\documentclass[12pt]{amsart}\usepackage{tikz-cd}
\usepackage{mathtools}
\usepackage{amsfonts}
\usepackage{graphicx}
\usepackage{verbatim}
\usepackage{textcomp}\usepackage{amssymb}
\usepackage{cite}
\usepackage{color}
\usepackage[all]{xy}
\usepackage{graphicx}
\usepackage{color}

\definecolor{NoteColor}{rgb}{1,0,0}

\newcommand{\note}[1]{\textcolor{NoteColor}{^\sharp 1}}

\usepackage{hyperref}
\usepackage{soul}
\hypersetup{
    colorlinks,
    citecolor=black,
    filecolor=black,
    linkcolor=black,
    urlcolor=black
}
\usepackage{tikz}
\usepackage[all]{xy}
\usepackage{graphicx}
\usetikzlibrary{backgrounds}
\newtheorem{theorem}{Theorem}[section]
\newtheorem{lemma}[theorem]{Lemma}
\newtheorem{proposition}[theorem]{Proposition}
\newtheorem{corollary}[theorem]{Corollary}   

\newtheorem{claim}[theorem]{Claim}
\allowdisplaybreaks
\newtheorem{notation}[theorem]{Notation}

\theoremstyle{definition}
\newtheorem{definition}[theorem]{Definition}
\newtheorem{remark}[theorem]{Remark}

\newtheorem{example}[theorem]{Example}
{{\sc Proof of Theorem~\ref{earthquakethm}}}
{{\sc q.e.d.} \\}

\newenvironment{thm:1relvw}
{{\sc Proof of Theorem~\ref{thm:1relvw}.}}
{{\sc q.e.d.} \\}
\newenvironment{genNTheorem 2}
{{\sc Proof of Theorem~\ref{genNTheorem 2}.}}
{{\sc q.e.d.} \\}
\newenvironment{proofconvvq}%
{{\sc Proof of Claim~\ref{convvq}.}}%
{{\sc q.e.d.} \\}
{{\sc Proof of Lemma~\ref{cd2general}.}}%
{{\sc q.e.d.} \\}
{{\sc Proof of Lemma~\ref{usordv}.}}%
{{\sc q.e.d.} \\}
{{\sc Proof of Theorem~\ref{thm:supptmeasure}.}}%
{{\sc q.e.d.} \\}
{{\sc Proof of Theorem~\ref{mainthmpushforw}.}}%
{{\sc q.e.d.} \\}
 \newenvironment{Length=totvar}%
{{\sc Proof of Theorem~\ref{Length=totvar}.}}%
{{\sc q.e.d.} \\}
 \newenvironment{convthmrepresentearth}%
{{\sc Proof of Theorem~\ref{convthmrepresentearth}.}}%
{{\sc q.e.d.} \\}
  \newenvironment{Length-totvarThur}%
{{\sc Proof of Theorem~\ref{Length-totvarThur}.}}%
{{\sc q.e.d.} \\}

  \newenvironment{thm:existence}%
{{\sc Proof of Theorem~\ref{thm:existence}.}}%
{{\sc q.e.d.} \\}
\numberwithin{equation}{section}

{{\sc Proof of Theorem~\ref{regularity}.}}%
{{\qed} \\}
\newenvironment{proof:main}%
{{\sc Proof of Theorem~\ref{ex}.}}%
{{\qed} \\}

\newenvironment{proof:pluriharmonic}%
    {{\sc Proof of Theorem~\ref{theorem:pluriharmonic}.}}%
  {{\qed} \\}  
  
  \newenvironment{proofof(iv)}%
    {{\sc Proof of $(iv)$.}}%
  {{\qed} \\}  
\newcommand{\norm}[1]{\left\Vert^\sharp 1\right\Vert}
\newcommand{\abs}[1]{\left\vert^\sharp 1\right\vert}
\newcommand{\R}{\mathbb R}

\newcommand{\HH}{\mathbb H}

\begin{document}
\title{Best Lipschitz maps and Earthquakes }
\author[Daskalopoulos]{Georgios Daskalopoulos}
\address{Brown Univeristy \\
Providence, RI}
\author[Uhlenbeck]{Karen Uhlenbeck}
\address{University of Texas Austin, TX and Institute for Advanced Study, Princeton, NJ}
\thanks{GD supported in part by NSF DMS-2105226}
\maketitle

\begin{abstract} 
This is the third paper (cf. \cite{daskal-uhlen1}, \cite{daskal-uhlen2}) in which we 
prove Thurston's conjectural duality between best Lipschitz maps and transverse measures (cf. \cite[Abstract]{thurston}).  
In \cite{daskal-uhlen2} we found a special class of best Lipschitz maps between hyperbolic surfaces  (infinity harmonic maps), which induce  dual Lie algebra valued transverse measures  with support on Thurston's canonical lamination. The present paper examines these Lie algebra valued  measures in greater detail. For any measured lamination we are led to define a  Lie algebra valued  measure and conversely  every Lie algebra valued transverse measure arrises from this process. Furthermore, we  show that such measures are  infinitesimal earthquakes. This construction provides a natural correspondence between best Lipschitz maps and earthquakes.
\end{abstract}

\section{Introduction}

The classical theorem of Kirszbraun-Valentine (1934)  states that a Lipschitz map from a subset $A$ of  $M_1=\R^n$ to   $M_2=\R^m$ can be extended to the entire space without increasing the Lipschitz constant (cf. \cite{kirszbraun}).  For $M_2=\R$ similar theorems  were obtained independently at about the same time by McShane and Whitney (cf. \cite{shane} and \cite{whitney}). The case of scalar valued functions has also been extensively studied  by analysts in connection to the infinity Laplacian (cf. \cite{aronson1}, \cite{aronson2}, \cite{arcrju}, \cite{crandal} and the references therein).

Our interest in the theory is based on the  1998 preprint by Bill Thurston \cite{thurston} (still unpublished), which constructs  a Teichm\"uller theory of hyperbolic surfaces coming from best Lipschitz maps $u: M_1=(M, g) \rightarrow M_2=(N, h)$ between hyperbolic surfaces homotopic to the identity (cf. \cite{thurston}).   The  construction is geometric: For a simple closed curve $\gamma$, define $l_g(\gamma)$ (resp. $l_h(\gamma))$  the length of the geodesic homotopic to $\gamma$ in $M $ (resp $N$).  Set 
\begin{equation}\label{defKintro}
K = \sup_\gamma \frac{l_h(\gamma)}{\l_g(\gamma)}.
\end{equation} 
The supremum in~(\ref{defKintro}) may not be attained at a simple closed curve but  at a measured lamination (not necessarily unique). However, there exists a {\it unique} geodesic lamination $\lambda$ (possibly without a measure of full support)  containing all such measured laminations. $\lambda$ is called {\it the canonical lamination} and it will play a crucial role in  this paper.

Given $\lambda$, Thurston constructs a {\it best Lipschitz map} $u$ homotopic to the identity with Lipschitz constant $L = K$ that stretches the leaves  of the lamination by a linear factor equal to $L$. 
  Everyone familiar with Teichm\"uller theory may realize the analogy with Teichm\"uller maps  stretching the leaves of quadratic differentials by a constant equal to the dilatation. 
  
  Gueritaud and Kassel in \cite{kassel} revisited and generalized Thurston's construction with a different approach: Given hyperbolic surfaces  $M$ and $N$, they first define the canonical lamination $\lambda$ as the intersection of the maximum stretch set of all best Lipschitz maps between $M$ and $N$ in a given homotopy class. They proceed to show the existence  of a best Lipschitz map that stretches $\lambda$ by $L$ and deduce the equality $K=L$ by approximation with closed geodesics \cite[Lemma 5.2 and Lemma 9.3]{kassel}.
  
The purpose of this paper is to study, in an explicit way, the connection between a special class of best Lipschitz maps (infinity harmonic) and measured  laminations.
  The main idea  is to exploit a duality principle between $L^\infty$-minimizers and measures. The study of dual functionals is fairly common in optimization theory and  is often referred to as ``Legendre-Fenchel duality"  or   ``max flow min cut principle" (cf. \cite[Abstract]{thurston}). For an  abstract theory of dual functionals in analysis see \cite[Chapter IV]{temam} or for the special case of conjugate harmonic functions  \cite {aronsonlin}. See also \cite{aidan2}. 
%


Before we describe in more detail our main results, we give a rough outline of \cite{daskal-uhlen1} and \cite{daskal-uhlen2}:
In our first paper, we studied the minima $u_p$ as $p \rightarrow \infty $ of the functional
\[
    J_p(f) =  \int_M |df|^p d\mu
\]
for $f:M_1 \rightarrow M_2$ where $M_1=M$ is a closed Riemannian manifold and $M_2=S^1$.  We obtained a limit $u_p \rightarrow u$ for a sequence $\{p \rightarrow \infty\}$ where $u$ is locally an infinity harmonic function.  Furthermore, after passing to a subsequence, the properly normalized closed one-forms
\[
               V_q =  k_p * |du_p|^{p-2}du_p \rightarrow V
\]
($*$ denotes Hodge star).  Here and throughout the rest of the paper $1/p+1/q=1$. The surprise in this paper was that
 $V$ is a transverse measure on the maximal stretch lamination   
\[
\lambda = \{x \in M; |du| = L \ \mbox{(:=the best Lipschitz constant)} \}.
\]

The second paper carries out this construction for $M_1 = M$, $M_2 = N$ both hyperbolic surfaces.  The best Lipschitz map is approximated by the minima  $u_p : M \rightarrow N$ of the functional
\[
               J_p(f) = \int_M  Tr Q(df)^p d\mu
\]
where $Q^2$ is the bilinear form on the pullback bundle $f^* TN$ constructed from $dfdf^T$ (Schatten-von Neumann norm).  We call $u_p$ {\it $p$-Schatten harmonic maps.} The  approximation $V_q$ to the dual measure predicted by Thurston is the closed Lie algebra valued 1-form
\[
               V_q = k_p*Q(u_p)^{p-2}du_p \times u_p
\]
where cross product is taken in the Lorenzian space $\R^{2,1}$.
As with the case $N = S^1$ (after passing to a subsequence),
$u_p \rightarrow u$
is a best Lipschitz map (which we call again {\it $\infty$-harmonic})
and      $V_q \rightarrow V.$
Here, the dual to the best Lipschitz map $u$, $V$ is a distribution with values in $T^*M$ tensor with the Lie algebra bundle of $N$. In addition, $V$ has finite  total variation, is supported on the canonical lamination and pushes forward by $u$ to a   measure on $N$. This is an example of what we call  a {\it Lie algebra valued transverse measure}.

 It is interesting to point out that, in the case of maps between surfaces,  the duality principle conjectured by Thurston naturally produces a Lie algebra valued measure (or an earthquake) and not a real valued transverse measure.  To produce an actual transverse measure,  requires further analysis. 

The following theorem is the easiest to state and provides an explicit relation between the best Lipschitz map, the canonical lamination and the transverse measure: 

\begin{theorem}\label{MMMainTHM}Let $M$, $N$ be closed hyperbolic surfaces and let $u:M \rightarrow N$ be an infinity harmonic map such that $u=\lim_{p \rightarrow \infty} u_p$  for a sequence  $u_p$ of $J_p$-minimizers in the same homotopy class. Then after passing to a subsequence  and after normalizing by positive constants $\kappa_p$, $|S_{p-1}| := Tr(Q(\kappa_p du_p)^p)$ converges weakly to a transverse measure 
$|S|$ on the canonical lamination $\lambda$  \footnote{in the topology literature (cf. \cite[Chapter 8]{thurston2}) it is customary to require that transverse measures are of full support; for us a transverse measure could be supported on a sublamination.} associated to the hyperbolic metrics on $M$, $N$ and the homotopy class of $u$.
\end{theorem}

Note that  $u$ and the transverse measure depend a-priori on the sequence $p \rightarrow \infty$ and we do not know if they are unique. In fact, we don't know uniqueness even for $N=S^1$. It is natural to conjecture that this process produces a unique measure. We are able to verify the conjecture only in the  special case when the canonical lamination consists of closed geodesics.

A  synopsis of the paper is as follows:  The present  section is an introduction, an outline of the paper and a statement of the main results.  Section~\ref{sect:affbundlinfinite} consists of the results we need from \cite{daskal-uhlen2} to make the above rough description of the results useful and rigorous: Here we review $p$-Schatten harmonic maps  and their limiting infinity harmonic map $u=\lim_{p \rightarrow \infty}u_p$. Next we define  the dual Lie algebra valued measures $W$ and $V$  produced by conservation laws from the hyperbolic symmetries of the domain and target. We also construct primitives $w$, $v$ to $W=dw$, $V=dv$ which we interpret as sections of bounded variation of affine vector bundles with linear structure  the Lie algebra bundle of the domain and the target respectively.

Section~\ref{prelimeart} contains  one of the  main results. 
We start by defining the notion of a transverse measure on a geodesic lamination with values in the Lie algebra of $SO(2,1)$. An example  is  the measure $W=dw$ on $M$ coming  from a best Lipschitz map $u$. 
 We also define {\it the standard  Lie algebra valued  transverse measure  associated to a measured geodesic lamination:}  if $(\lambda, \mu)$
 is a measured geodesic lamination, we set $dw=Bd\mu$ where $B$ is the generator of the geodesic flow of $\lambda$.  This expression is independent of the local orientation of the lamination. We  prove  that any Lie algebra valued  transverse measure on a geodesic lamination arises from a (real valued) transverse measure by the above process:

\begin{theorem}\label{thmrepresentearth} Let $M$ be a closed hyperbolic surface and  $W=dw$  a Lie algebra valued  transverse measure with support in a  geodesic lamination $\lambda$ of $M$. Then  there exists a unique transverse measure $\mu$ with support in  $\lambda$  such that  $dw=Bd\mu$ where B is the generator of the geodesic flow of $\lambda$. Furthermore, 
$\mu=*(\omega_{mc} \wedge dw)^\sharp.$ Here, $\omega_{mc}=dx \times x$ denotes the Mauer-Cartan form on $M.$ (Note that multiplying a measure by a smooth form makes sense.)
  \end{theorem}

 We end Section~\ref{prelimeart} by proving the following theorem  which relates  mass (total variation) to length:
 \begin{theorem}\label{Length=totvar}Let $(\lambda, \mu)$ be a measured geodesic  lamination on the hyperbolic surface $(M,g)$ and $dw=Bd\mu$ a Lie algebra valued transverse measure. The mass (total variation) $||dw||_{meas}$ satisfies
\[
||dw||_{meas}=  2 l_g(\lambda).
\]
\end{theorem}

In  Section~\ref{sect:KL} we apply the results of the previous sections to best Lipschitz maps and the Lie algebra valued transverse measures $dw$, $dv$ on $M$, $N$ associated to an infinity harmonic map $u$. First we show for $dw$:  

\begin{theorem}\label{mainthmnotdom} Let $M$, $N$, $\lambda$  and $u$ be as before  and  $W=dw$  the  measure on $M$ coming from the Noether current of $u$.  Then  
 $dw$ is a Lie algebra valued transverse measure on $\lambda$ and $dw=B d\mu$ where $\mu=2|S| :=2\lim_{p \rightarrow \infty} Tr(Q(\kappa_p du_p)^p)$.  
 \end{theorem}
  
 For $dv$ we show:
 \begin{theorem}\label{thm:1relvw} Let $M$,$N$, $u$, $\lambda$  as in the previous theorem. Then
  \begin{equation}
dv=B^\land(u) d\mu
\end{equation}
where $B^\land$ is the generator of the geodesic flow of the geodesic lamination $\lambda^{\land}=u(\lambda)$ and $\mu=2|S|=*(\omega_{mc}, dw)^\sharp$. 
 \end{theorem}
%
%
%
In order to prove Theorem~\ref{thm:1relvw} we push $dv$ forward by the best Lipschitz map $u$ to  obtain a  Lie algebra valued  transverse measure $dv^{\land}$ on $\lambda^\land$ and  apply Theorem~\ref{thmrepresentearth}. 
Combining Theorem~\ref{mainthmnotdom} with Theorem~\ref{thm:1relvw} we conclude that $dw=B d\mu$ and $dv=B^\land(u) d\mu$ are both rank one Lie algebra valued measures which are {\it multiples of the same real valued measure transverse} $\mu$. The Lie algebra valued functions $B$ and $B^\land(u)$ are Lipschitz and rank one. This is  not a surprise in view of Alberti's rank one theorem for singular vector valued measures (cf. \cite{alberti}).

In Section~\ref{lamearth} we  describe  the connection between Lie algebra valued transverse measures and earthquakes. We show that the  Lie algebra valued transverse measure $dw=Bd\mu$  defines  an infinitesimal earthquake. 
 Earthquakes are changes in the hyperbolic structure constructed by cutting along  geodesics, twisting and regluing. As the twist goes to zero, we obtain  tangent vectors to Teichm\"uller space which are called  {\it infinitesimal earthquakes}. Equivalently,  infinitesimal earthquakes can be viewed as  cohomology classes in the Lie algebra cohomology $H^1(\pi_1(M), {\mathfrak g}_\sigma)$ representing tangent vectors to the space of conjugacy classes of discrete, faithful $SO^+(2,1)$ representations:

 \begin{theorem}\label{convthmrepresentearth} {\it Let $\sigma: \pi_1(M) \rightarrow SO^+(2,1)$ be a representation  defining a hyperbolic structure   and $(\lambda, \mu)$ a measured geodesic lamination on a closed surface $M$. Assume $dw=Bd\mu$ is a Lie algebra valued  transverse measure on  $\lambda$. Then the cohomology class of $dw$ in $H^1(\pi_1(M), {\mathfrak g}_\sigma)$ is equal to the infinitesimal earthquake of  $(\lambda, \mu)$ at $\sigma$.}
 \end{theorem}

 Our proof is based on the basic fact that  infinitesimal earthquakes  are dual to the derivative of length. We prove the theorem by first considering  the case when the lamination consists only of closed geodesics. Using results of Kerckhoff, we then pass to the limit to treat general laminations. In a subsequent paper, we will show that earthquake maps (cf. \cite{thurston3}) have a formulation as Lie Group valued transverse measures, and show that  $dw$ is geometrically the infinitesimal deformation of an earthquake map.

 In this paper, we consider both domain and target spaces to be hyperbolic surfaces.  There is considerable interest in allowing more general targets.  Our results go over with very little change to the case of $SO^+(n,1)/SO(n)$ and the Lipschitz constant $L > 1$.  For $L < 1$, the analysis carries over entirely. However, we do not know what replaces the measured lamination $(\lambda, \mu)$; it is quite possible the analysis of $dv$ and $dw$, which are well-defined, will be helpful in determining this.  For $SL(n,R)/SO(n)$, one can, of course, define infinity harmonic maps. Noether's theorem still applies, but the  support argument is harder to execute.  Best Lipschitz may not necessarily be the correct approach. The goal would be to define an object similar to  the map $u$ and the measures $dv$ and $dw$.

{\bf Acknowledgements.}
We would like to thank Aidan Backus for pointing out Anzellotti's paper  and Mike Wolf for sending us the preprint  \cite{wolf}.  We would also like to thank Jeff Danciger, Fanny Kassel and  Athanase Papadopoulos for useful discussions during the preparation of this manuscript.

\section{Preliminaries}\label{sect:affbundlinfinite} 
This section contains some background material on best Lipschitz  and infinity harmonic maps  and their dual measures. For further details we refer the reader to \cite{daskal-uhlen2}. We also review some basic facts about 1-cocycles, affine actions and the De Rham isomorphism on Lie algebra cohomology.  This allows us to interpret  the dual measure  to an infinity harmonic map as a derivative of a global section  of bounded variation of an affine vector bundle.
The measure, viewed as a closed 1-current, represents a well defined cohomology class in the first Lie algebra cohomology group.


\subsection{The $p$-approximation of infinity harmonic maps} \label{paprox} {\it In this section, as well as in the rest of the paper,  $(M,g)$ and $(N,h)$ are closed hyperbolic surfaces.} For   $2< p < \infty$ and $f \in W^{1,p}(M,N) $ we minimize the $p$-Schatten norm of the derivative of $f$ in a fixed homotopy class. More precisely, let $Q(df)^2 :=dfdf^T \in End(f^*(TN))$ and let $u_p$ be the minimizer of
\[
 J_p(f)=||df||^p_{sv^p}=\int_M |df|_{sv^p}^p*1=\int_M TrQ(df)^p*1.
\]
 The minimizer exists, is unique (unless it maps onto a single geodesic) and satisfies the Euler-Lagrange equations
\begin{equation*}\label{eqn:firstvar3}
D^*Q(du_p)^{p-2}du_p=0
\end{equation*}
where $D=D_{u_p}$ is the pullback of the Levi-Civita connection on $u_p^*(TN) $ \cite[Proposition 2.6 and Corollary 2.11]{daskal-uhlen2}.

The important property of the $u_p$'s is that, after passing to a subsequence, they limit to a {\it best Lipschitz map} i.e a map that minimizes the global Lipschitz constant in the homotopy class: 
 \begin{theorem}\label{thm:existence}\cite[Proposition 6.5 and Lemma 6.9]{daskal-uhlen2}   Given  a sequence  $p \rightarrow \infty$, there exists a subsequence  (denoted also by $p $) and a sequence of $J_p$-minimizers 
$u_p: M \rightarrow N$   in the fixed homotopy class and $u_p \rightarrow u$ in $C^0$ and weakly in $W^{1,s}$ for all $s$. Furthermore,  $u$ is a best Lipschitz map in the homotopy class and
\begin{equation}\label{pintconvto}
\lim_{p \rightarrow \infty} J_p^{1/p}(u_p) =  L
\end{equation}
where $L$ is the best Lipschitz constant. We call $u$ an {\it infinity harmonic map}.
\end{theorem}

\begin{remark} Note that in this paper  infinity harmonic maps are limits as $p \rightarrow \infty$ of $p$-Schatten harmonic maps and not, as some people call them,  limits of minimizers of the $L^p$-norm. The latter do not converge to a best Lipschitz map and the geometric meaning is unclear.
\end{remark} 

 Of crucial importance in our study is  the 
{\it canonical lamination} $\lambda=\lambda(g,h)$ associated to the hyperbolic metrics and the homotopy class defined by Thurston and Gueriteaud-Kassel:
Given a best Lip map $f: (M,g) \rightarrow (N,h)$ in a given homotopy class, let 
\[
 \lambda_f=\{ x \in M: |df(x)| =L \} 
\]
 be its maximum stretch locus. Let $\mathcal F$ denote the collection of all best Lipschitz maps in the homotopy class and set
\[
\lambda=\cap_{f \in \mathcal F} \lambda_f.
\]
Then $\lambda$ is a geodesic lamination \cite[Lemma 5.2]{kassel} which, in the case of the homotopy class of the identity,  is equal to Thurston's chain recurrent lamination  associated to the hyperbolic structures $g,h$. (cf. \cite[Theorem 8.2]{thurston} and \cite[Lemma 9.3]{kassel}).

\subsection{The hyperboloid in $\R^{2,1}$}\label{hyperboloid} 
The hyperbolic plane $\HH$ viewed as the unit hyperboloid in $\R^{2,1}$. Denote the indefinite metric on  $\R^{2,1}$ inducing the hyperbolic metric on $\HH$ by $(,)^\sharp$.
We review some basic geometry.   Let $e^\sharp  = diag(1,...,1,-1)$.  For $X \in \R^{n,1}$ let    $X^\sharp  = (e^\sharp X)^T$. The inner product in $\R^{n,1}$ is
\[  
( X,Y)^\sharp  =  X^\sharp Y =Y^\sharp X.
\] 
The associated transpose on linear maps $B$  is $B^\sharp  = e^\sharp B^Te^\sharp$.  The group $SO (n,1)$ consists of $(n+1) \times (n+1)$ matrices with determinant one that preserve the inner product. Equivalently, $g \in G$ if $det g=1$ and  $g^{-1} = g^\sharp$.  $B \in \mathfrak s \mathfrak o(n,1)$ iff $Tr B=0$ and  $B^\sharp  = -B$.   For $X,Y \in \R^{n,1}$,  $XY^\sharp$  is a matrix and $YX^\sharp  - X^\sharp Y$ is a skew symmetric matrix (with respect to $^\sharp$). Define
\begin{equation}\label{wedgedef}
X \times Y=YX^\sharp  - X Y^\sharp \in  \mathfrak g=so(n,1).
\end{equation}
  Let  $(A,B)^\sharp= Tr A B $ denote the Killing form on $\mathfrak g$. 
 With respect to a Cartan decomposition 
 $ \mathfrak g=so(n,1)$,
  write $A \in \mathfrak g$ as 
 \begin{equation*}
A = 
\begin{pmatrix}
W & X  \\
X^T & 0
\end{pmatrix}
\end{equation*} 
 for $X$ arbitrary vector and $W=-W^T$.
  It follows
 $(A,A)^\sharp= Tr WW+2|X|^2 $
  and the metric has signature $(n, n(n-1)/2)$. In particular for $n=2$, the Killing form is (up to a constant) a flat Lorenzian metric on $\mathfrak g=so(2,1)$ of signature $(2,1)$  isometric to $\R^{2,1}$.
 Let
\[
\HH^n = \HH = \{X \in \R^{n,1}: ( X,X )^\sharp  = -1 \ and \ X_{n+1} \geq 1 \} 
\] 
and $G=SO^+(n,1)$ the index two subgroup of $SO(n,1)$ preserving $\HH$.

%
For $X$ a point in $\HH$, let $\Pi(X)$ be the orthogonal in $(, )^\sharp$  projection onto $T_X\HH$ and $\Pi^\perp(X)$ be the orthogonal projection onto $X$.  
Then
\begin{equation*}
 I = \Pi(X)+ \Pi(X)^\perp,  \ \Pi^\perp(X)=-XX^\sharp 
 \end{equation*}
and
\[
 \Pi(X)=I+XX^\sharp.
\]
For $v \in \R^{1,n}$,  denote the projection
\begin{equation} \label{proj}
\Pi(X)v=v_X=v+(v,X)^\sharp X.
 \end{equation}
This is similar to the formulas for $S^n$ in $\R^{n+1}$ with the change in sign in $XX^\sharp$ due to the indefinite metric in $\R^{n,1}$.   
  \begin{lemma}\label{metricHemb}The inner product $(,)^\sharp$ restrictred to the tangent bundle of $\HH$ is a Riemannian metric which agrees with the standard metric of the hyperbolic space.
  \end{lemma}

\subsection{Conservation laws}\label{conservlawsrarg} In \cite{daskal-uhlen2} we described  conservation laws associated to the symmetries of the domain and the target. Let $G=SO^+(2,1)$ acting on the hyperbolic plane $\HH$ viewed as the unit hyperboloid in $\R^{2,1}$. Denote the indefinite metric on  $\R^{2,1}$ inducing the hyperbolic metric on $\HH$ by $(,)^\sharp$. Identify the universal covers of $M$ and $N$ with $\HH$ and assume  the hyperbolic metrics $g$ and $h$ are defined by  discrete and faithful representations 
\[
\sigma: \pi_1(M) \rightarrow G=SO^+(2,1)\ \mbox{and} \ \rho: \pi_1(N) \rightarrow G=SO^+(2,1)  
\]
respectively.

An element $w$ of the Lie algebra $\mathfrak g=so(2,1)$  defines a Killing vector field on $\HH$ by setting   $w(X)=w X \in T_X\HH$. Here $X \in \HH \subset \R^{2,1}$ and  $w$ acts on $X$ as a skew-adjoint endomorphism with respect to $(,)^\sharp$. 
This defines a map
\[
\alpha_X: \mathfrak g \rightarrow T_X\HH \subset \R^{2,1}, \ \  w \mapsto wX.
\]
There is a right inverse 
\[
\beta_X: T_X\HH \rightarrow \mathfrak g
\]
 of $\alpha_X$  which is best described via {\it the cross product}.
 
For $X \in \R^{2,1}$ set    $X^\sharp  = (e^\sharp X)^T$
where $e^\sharp  = diag(1,1,-1)$. The inner product in $\R^{2,1}$ is 
$( X,Y)^\sharp  =  X^\sharp Y =Y^\sharp X$
 and the associated transpose on linear maps $B$  is $B^\sharp  = e^\sharp B^Te^\sharp$.  Then $g \in G =SO^+(2,1)$ preserves the inner product  and satisfies $g^* = g^\sharp$.  $B \in \mathfrak g=so(2,1)$ satisfies $B^\sharp  = -B$.   Note that, for $X,Y \in \mathfrak g$,  $YX^\sharp  - X^\sharp Y$ is a skew symmetric matrix (with respect to $^\sharp$). Define
\begin{equation}\label{wedgedef}
X \times Y=YX^\sharp  - X Y^\sharp \in  \mathfrak g=so(2,1)
\end{equation}
and set
\[
\beta_X(v)= v \times X.
\]
Here both $v$ and $X$ are considered as points of $\R^{2,1}$. We can now let $X$ vary over $\HH \subset \R^{2,1}$. 

By letting $X$ vary and keeping track of the equivariance,  the above pointwise construction
  yields global maps 
\[
\alpha: ad(E_\sigma) \rightarrow TM \ \mbox{and} \ \beta: TM \rightarrow ad(E_\sigma)
\] 
between the tangent bundle of $M=\HH/\sigma(\pi_1(M))$ and the  Lie algebra bundle $ad(E_\sigma)=\HH \times_{Ad(\sigma)} \mathfrak g$.

We can generalize this construction for  maps between the surfaces as follows. Let $f: M=\HH/\sigma(\pi_1(M)) \rightarrow N=\HH/\rho(\pi_1(M))$ be a Lipschitz map and let
$\tilde f: \tilde M =\HH^2 \rightarrow \tilde N =\HH^2$ an equivariant  lift of $M$ intertwining the representations $\sigma$ and $\rho$.
 Define 
 the flat $\R^{2,1}$-bundle $E_\rho=\tilde M \times_\rho \R^{2,1}$ and the flat Lie algebra bundle $ad(E_\rho)=\tilde M \times_{Ad(\rho)} \mathfrak g$ over $M$.   
Letting $X=f(x)$ where $x \in \tilde M= \HH$, the pointwise construction described earlier gives  global maps 
\begin{equation}\label{defn:alphabetaf}
\alpha_f: ad(E_\rho) \rightarrow f^*(TN) \ \mbox{and} \ \beta_f: f^*(TN) \rightarrow ad(E_\rho).
\end{equation}
Note that for $\rho=\sigma$ and $f=id: M \rightarrow M$, we recover $\alpha_{id}=\alpha$ and $\beta_{id}=\beta$.
The other case of interest is when  $f=u: M \rightarrow N$ is an infinity harmonic map or $f=u_p$ its $p$-approximation. For further details we refer to  \cite[Section 3.3]{daskal-uhlen2}.

For $u_p$ a minimizer of $J_p$, let  
\begin{equation}\label{rec1}
S_{p-1}=Q(du_p)^{p-2}du_p, \  \ \ V_q=*\beta_{u_p}S_{p-1}=*S_{p-1} \times u_p
\end{equation}
Here, as in the entire paper, 
\begin{equation}\label{form:conjugate}
1/p+1/q=1.
\end{equation}
 The conservation law coming from the local symmetries of the target is:

\begin{proposition}\label{Prop:Elag-bund} \cite[Theorem 3.9]{daskal-uhlen2}  The Noether current $V_q$ is closed with respect to the flat connection $d$ on $ad(E_\rho)$, i.e. $dV_q=0$.   
 \end{proposition}

There is a corresponding conservation law associated to the bundle $ad(E_\sigma)$ of the local symmetries of the domain $M$. Set 
\begin{equation}\label{rec2}
|S_{p-1}|=Tr Q(du_p)^p.
\end{equation} 
The energy-momentum tensor is
\begin{equation}\label{rec3}
T_q=T_q(du_p) =  (S_{p-1} \otimes du_p)^\sharp - 1/p|S_{p-1}| g.
\end{equation}
The notation $(.  \otimes   .)^\sharp$ means that we are taking tensor product in the domain and contracting by the metric on the target.
Let
\begin{equation}\label{rec4}
W_q=*\beta(T_q)=*T_q \times id 
\end{equation}
a section of $ad(E_\sigma)$.
Then:

\begin{proposition}\label{Prop:NOTD} \cite[Theorem 4.4 and Proposition 4.8]{daskal-uhlen2}  The energy momentum tensor $T_q$ is divergence free with respect to the Levi-Civita connection. Equivalently,  $W_q$ is closed with respect to the flat connection $d$ on $ad(E_\sigma)$, i.e. $dW_q=0$.   
 \end{proposition}
 
 The following proposition is immediate from the definitions:
 \begin{proposition}\label{prop:relvw} The tensors $T_q$ and $V_q$ are related by
\[
-2T_q=(*V_q \otimes du_p \times u_p)^\sharp.
\]
\end{proposition}

We conclude with the following pointwise construction and its global version that will be used throughout the paper:
For $X \in \HH$ and $v \in \R^{2,1}$ denote by  $v_X$ the orthogonal projection of $v$ onto the tangent space $T_X \HH$ with respect to the indefinite metric on $\R^{2,1}$.
Similarly, for $v \in \mathfrak g$ denote by  $v_X$ the orthogonal projection of $v$ onto the positive definite part of $\mathfrak g$ with respect to the indefinite metric on $\R^{2,1}$.
For $\xi \in \Omega^1(E_\rho)$, let $\xi_f$ denote pointwise the orthogonal projection onto the subbundle $f^*(TN) \subset E_\rho$. More precisely, let 

Similarly, for $\phi \in \Omega^1(ad(E_\rho))$, let $\phi_f$ denote pointwise the orthogonal projection (using  the metric $(,)^\sharp$) of $\phi$ onto the image of $\beta_f$.
 
\subsection{The limiting measures}\label{limleassecti} 

As $p \rightarrow \infty$ and the minimizers $u_p$ approach a best Lipschitz map $u$,   the tensors defined in (\ref{rec1}),
(\ref{rec2}), (\ref{rec3}) and (\ref{rec4}),  when appropriately normalized, converge weakly to  Radon measures. The normalizations are defined in \cite[Section 7.1]{daskal-uhlen2}. For the purpose of this paper we remind the reader  that we choose constants $\kappa_p \rightarrow L^{-1}$ such that the $p$-Schatten norm of $U_p=\kappa_pdu_p$ is one. Here $L$ is the best Lipschitz constant. We define the renormalized tensors $S_{p-1}$, $V_q$, $T_q$ and $W_q$ using  $U_p$ instead of $du_p$. Then

 \begin{theorem}\label{thm:limmeasures0} \cite[Theorems 7.4, 7.6 and 7.8]{daskal-uhlen2}  Given a sequence $p \rightarrow \infty$ ($q \rightarrow 1$), there exists a subsequence (denoted again by $p$, $q$), 
a real-valued positive Radon measure $|S|$    
and Radon measures $S$, $V$, $T$ and $W$ with values in $T^*M \otimes E_\rho$, $T^*M \otimes ad(E_\rho)$, $T^*M \otimes E_\sigma$ and  $T^*M \otimes ad(E_\sigma)$ respectively such that, after rescaling:
\begin{itemize}
\item $(i)$ $ |S_{p-1}| \rightharpoonup |S|$, $ S_{p-1} \rightharpoonup S$, $ V_q \rightharpoonup V$, $ T_q \rightharpoonup T$ and $ W_q \rightharpoonup W$
 \item $(ii)$$dV=0$,  $dW=0$ with respect to the flat connection on $ad(E_\rho)$ and $ad(E_\sigma)$ respectively
 \item $(iii)$ The supports of $S$ and $V$ are equal and contained in the support of $|S|$. 
\item $(iv)$ The supports of $T$, $W$ and $|S|$ are all equal. Furthermore, $*(\omega_{mc} \wedge W)^\sharp=2|S|.$
 \item $(v)$ The support of $|S|$ (and thus of all the measures) is contained in the canonical geodesic lamination $ \lambda$ associated to the hyperbolic metrics $g$, $h$ and the homotopy class
 \item $(vi)$ $S_u=S$, $V_u=V$, $T_{id}=T$, $W_{id}=W$
 \item $(vii)$  The mass (total variation) of $|S|$, $S$ and $T$ is one and of $V$, $W$ is two.
 \end{itemize}
\end{theorem} 

Some  explanation for the notation is needed for $(iv)$, $(vi)$ and $(vii)$. First, we explain the Mauer-Cartan form $\omega_{mc}$. The cross product defines a map $\times: T\HH  \rightarrow \mathfrak g$ which intertwines the action of $SO^+(2,1)$ on $\HH$ with the adjoint action on $\mathfrak g$. Therefore cross product defines a section of  $T^*\HH \otimes \mathfrak g$ which by equivariance descends to $\omega_{mc} \in \Omega^1(M, ad(E_\sigma)).$
 Sometimes we use the shorthand $\omega_{mc}=dx \times x$. 

\begin{definition}\label{projcurr} Let $f: M \rightarrow N$ be  a Lipschitz map, $\gamma$ a distribution  with values in $T^*M \otimes E_\tau$ and $\zeta$ a distribution  with values in  $T^*M \otimes ad(E_\tau)$. We define $\gamma_f$ and $\zeta_f$ by setting
\[
\gamma_f(\xi)=\gamma(\xi_f) \ \ \ \zeta_f(\phi)=\zeta(\phi_f)
\]
for all  $\xi \in \Omega^1(E_\tau)$ and $\phi \in \Omega^1(ad(E_\tau)).$
\end{definition}

The condition  $\gamma_f=\gamma$ means that $\gamma$ is a measure with values on the subbundle $T^*M \otimes f^*TN \subset T^*M \otimes  E_\tau.$ The condition  $\zeta_f=\zeta$  means that $\zeta$ has values on the $T^*M \otimes  Im(\beta_f) \subset ad(E_\tau).$ Recall that  $Im(\beta_f)$ is the subbundle of the Lie algebra bundle where the metric is positive definite.


Finally we discuss $(vi)$ and the term mass. 
 For the real valued measure $|S|$, mass is the same as the usual total variation. However, for the vector valued measures the definition of mass is more sensitive to which norm we are using. In our case the right definition  on test functions is to use the operator norm ($\infty$-Schatten norm) and not the usual $L^\infty $ norm. 
 
\begin{definition}\label{measrcurr} Let $f$, $\gamma$ and $\zeta$ as in the previous definition.  We define the {\it mass}  of $\gamma$ and $\zeta$ with respect of $f$ 
\[
||\gamma||_{meas,f} =\sup \{ \gamma(\xi): \xi \in \Omega^1(E_\tau), \xi=\xi_u, s( \xi) \leq 1 \}.
\]
and 
\begin{eqnarray*}
||\zeta||_{meas,f} &=& \sup \{ \zeta(\xi \times f): \xi \in \Omega^1(ad(E_\tau)),  s( \xi_f) \leq 1 \}\\
&=& \sup \{\zeta(\phi): \phi \in \Omega^1(ad(E)), \phi=\phi_f, s( \phi) \leq \sqrt2\}.
\end{eqnarray*}
(In the above $s(\xi)$ denotes pointwise the largest singular value of $\xi$ viewed as homomorphism from $TM$ to $E$.  Similarly for $s(\phi)$ (cf. \cite[Section 3]{daskal-uhlen2}).  The factor of $\sqrt2$ comes from the fact that  $s(\xi \times f)=\sqrt2 s( \xi_f)$ (see again 
\cite[Section 3]{daskal-uhlen2}).
If $f=id$, we denote $||S||_{meas,f}=||S||_{meas}$ and $||V||_{meas,f}=||V||_{meas}$.
\end{definition}

\begin{remark} Later we will improve Theorem~\ref{thm:limmeasures0}$(iii)$ and show that the supports of $S$, $V$ are  equal to the support of $|S|$.  
\end{remark}

%

\subsection{1-cocycles and the De Rham isomorphism}\label{cocDer} 

 Let $G=SO^+(2,1)$ with Lie algebra $\mathfrak g=so(2,1)$ and consider the adjoint representation
$Ad: G \rightarrow Ad(G) \subset GL(\mathfrak g).$
 We also fix a  representation $\rho: \pi_1(M) \rightarrow G$ and we let $ad(E_\rho)$ the flat Lie algebra bundle as before. Let $H^1(M, ad(E_\rho))$ denote the De Rham cohomology group with values in $ad(E_\rho)$ associated to the flat connection $d=d_\rho$ on $ad(E_\rho)$. 

A 1-cocycle $\alpha \in Z^1(\pi_1(M), \mathfrak g_{Ad(\rho)})$ is defined by the {\it cocycle condition}
\begin{equation}\label{cohoalpha}
\alpha(\gamma_1\gamma_2)=Ad(\rho(\gamma_1))\alpha(\gamma_2)+\alpha(\gamma_1).
\end{equation}
A 1-coboundary  $\alpha \in B^1(\pi_1(M), \mathfrak g_{Ad(\rho)})$ is a 1-cocycle of the form
$\alpha=\alpha^0-Ad(\rho(\gamma))\alpha^0$,
where $\alpha^0$ is a fixed element in the Lie algebra. We let 
\[
H^1(\pi_1(M), \mathfrak g_{Ad(\rho)})=\frac{Z^1(\pi_1(M), \mathfrak g_{Ad(\rho)})}{B^1(\pi_1(M), \mathfrak g_{Ad(\rho)})}
\]
 and represent the cohomology class of the cocycle $\alpha$ by
$[\alpha] \in H^1(\pi_1(M), \mathfrak g_{Ad(\rho)}).$
A 1-cocycle defines an action of $\pi_1(M)$ by isometries on  $\mathfrak g$ given by
\begin{eqnarray}\label{affACt}
P \in \mathfrak g \mapsto Ad(\rho(\gamma))P+ \alpha(\gamma) \in \mathfrak g.
\end{eqnarray}
Such an action is called an {\it affine representation  with linear part  $Ad(\rho)$ and affine part $\alpha$} (cf. \cite{dacing}).

Given a cocycle $\alpha$, there exists a smooth map $\tilde f: \tilde M \rightarrow \mathfrak g$ 
such that
\begin{equation}\label{aqeqiov8}
\tilde f(\gamma z)=Ad(\rho(\gamma)) \tilde f( z)+\alpha(\gamma). 
\end{equation}
To construct such  $\tilde f$ we 
first pick $\alpha^0 \in \mathfrak g$, a fundamental domain $F$, define $\tilde f_0(F)=\alpha^0$ and extend equivariantly. This map is equivariant but not continuous, so we need to apply a mollifier. The mollification process in \cite[Proposition 2.6.1]{korevaar-schoen1} yields a Lipschitz equivariant map $\tilde f$. 
If we want a smooth map, we can just take the harmonic map homotopic to $\tilde f$.
Note that $\tilde \phi=d \tilde f$ is $Ad(\rho)$-equivariant, hence it descends to a closed 1-form with values in $ad(E_\rho)$ which we will denote by $df \in \Omega^1(M, ad(E_\rho))$. Note that this is a slight abuse of notation as $df$ is not exact. We now have the following analogue of the De Rham isomorphism for Lie algebra cohomology:

\begin{theorem}\label{cohoisp}The map
\[
 H^1(\pi_1(M), \mathfrak g_{Ad(\rho)})\rightarrow H^1(M, ad(E_\rho)), \ [\alpha] \mapsto [df]
\]
is an isomorphism. 
\end{theorem}

\begin{proof}
We will briefly sketch a proof. First note, that by adding a coboundary $\alpha^0-Ad(\rho(\gamma))\alpha^0$ to $\alpha$ changes $\tilde f$ to $\tilde f+ \alpha^0$ and thus  $df$ remains unchanged. The map is thus well-defined. To show that the map is onto, let $\phi \in \Omega^1(M, ad(E_\rho))$ closed with respect to the flat connection. Lifting to the universal cover, we can write $\tilde \phi=d \tilde f$ for $\tilde f: \tilde M \rightarrow \mathfrak g$. 
%
Since $\gamma^*(d \tilde f)=Ad(\rho(\gamma))(d \tilde f)$, the map 
 $ \tilde f$ satisfies  $ \tilde f(\gamma x)=Ad(\rho(\gamma)) \tilde f(x)+ \alpha(\gamma)$. This equivariance implies  cocycle condition (\ref{cohoalpha}) and proves that $[\phi]$ is the image of $\alpha$. Thus, the map is onto. In order to show the map is 1-1 assume $[\alpha]$ maps to $[\phi]=0$. Then $\phi=df$ for a section $f$ of $ad(E_\rho)$. Lifting to the universal cover, 
$\tilde f(\gamma x)=Ad(\rho(\gamma)) \tilde f(x)$ and thus $\alpha$ is cohomologous  to the zero cocycle.  
\end{proof}

\subsection{Constructing global primitives of $V$ and $W$}\label{primV}
We have already seen that  the 1-currents $V$ and $W$ are $d$-closed, with respect to the flat connection $d$ on $ad(E_\rho)$ and  $ad(E_\sigma)$ respectively (cf. Theorem~\ref{thm:limmeasures0}).  Thus, there exist local primitives $v$, $w$ of $V$, $W$, i.e local Lie algebra valued functions of bounded variation such that $V=dv$ and $W=dw$ (cf. \cite[Section 7]{daskal-uhlen2}).
 In this section we study global properties of  $v$ and $w$.
 
In \cite{daskal-uhlen1} we encountered a similar situation for the linear problem. There, $V$ was a closed (real  valued)  1-current on $M$ of finite total variation. Pulling back to the universal cover, we wrote $\tilde V =d \tilde v$ where $\tilde v$ is a function of locally bounded variation satisfying the equivariance condition
\begin{equation}\label{equivline}
\tilde v(\gamma x)=\tilde v(x)+ \alpha(\gamma) \ \ \ \gamma \in \pi_1(M), x \in \tilde M.
\end{equation}
 We interpreted this by saying that $\tilde v$ satisfies the affine condition for some representation $\alpha: \pi_1(M) \rightarrow \R$, or equivalently that $\tilde v$ induces a section $v$ of an affine line bundle on $M$.
 The situation we have here is very similar, but with the additional complication that $V$ has values in $ad(E_\rho)$ instead of the trivial line bundle.  This affects the equivariance condition of the primitive  $\tilde v$ and (\ref{equivline}) becomes (\ref{aqeqiov8})  where $\alpha$ is a 1-cocycle. 
 
 Consider  $V=V_q$ for $q>1$ first. As in the proof of Theorem~\ref{cohoisp}, we can lift to the universal cover to find $\tilde v_q$ such that $dv_q=V_q$ in the distribution sense. Here $\tilde v_q$ satisfies the equivariance condition~(\ref{aqeqiov8})  for a  1-cocycle $\alpha_q$.  Note that one needs to justify the lack of regularity as $V_q \in L^s$ for all $s$ but is not necessarily smooth (cf. \cite[Corollary 5.18]{daskal-uhlen2}).
We can easily overcome this problem by smoothing $V_q$ in the same cohomology class by running the heat equation and applying Theorem~\ref{cohoisp}. This construction is fairly standard so we skip the details.
\begin{theorem} \label{thm:limmeasures04r} Consider a sequence $\tilde V_q=d \tilde v_q$ where  $q \rightarrow 1$ and $\tilde v_q$ satisfies condition (\ref{aqeqiov8})  for a 1-cocycle $\alpha_q$. Assume $V_q$ has $q$-Schatten norm one.
Then there exists a subsequence (denoted again by $q$) such that 
 \begin{itemize}
 \item $(i)$ There exists a map of locally bounded variation $\tilde v: \tilde M \rightarrow \mathfrak g$ satisfying 
 condition (\ref{aqeqiov8})  for a 1-cocycle $\alpha$
such that $\tilde v_q$ converges to $\tilde v$ {\it weakly} in $BV_{loc}$ and {\it strongly} in $L^s_{loc}$ for  $s <2$. 
 \item $(ii)$ The cocycles converge, i.e $ \alpha_q \rightarrow \alpha.$
  \item $(iii)$   $d\tilde v$ descends to a measure on $T^*M \otimes ad(E_\rho)$  with support in the canonical lamination $\lambda$. We write  $dv=V$ and its
cohomology class $[dv]=[\alpha] \in H^1(M, ad(E_\rho))$. 
 \end{itemize}
  \end{theorem}

\begin{remark}
 The assumption that  $V_q$ has $q$-Schatten norm one  follows  the normalization assumption in Theorem~\ref{thm:limmeasures0}.
 (cf. \cite[Formula (7.4)]{daskal-uhlen2}).
\end{remark}

  \begin{proof}   The proof  follows in the spirit of  \cite[Theorems 4.1 and 4.3]{daskal-uhlen1} with the additional complication that the forms here have values in the adjoint representation. First, by adjusting the maps $\tilde v_q$ by an additive constant in the Lie algebra, we may assume that the average values of $\tilde v_q$ in a fundamental domain $F$ are uniformly bounded. As in the proof of Theorem~\ref{cohoisp} changing $\tilde v_q$ by a constant has the effect of adjusting each $\alpha_q$ by a coboundary and this has no effect in cohomology. 

By the Poincare inequality, $\tilde v_q$ are  uniformly bounded in $W^{1,1}(F)$. Hence there exists  $\tilde v \in BV(F)$ such that
$\tilde v_q \xrightharpoonup{BV(F)} \tilde v.$
It also converges strongly in $L^s(F)$ for  $s <2$ by Sobolev embedding. We next claim:

 \begin{claim}\label{convvq} Assume that  $\tilde v_q$ as before converge to $\tilde v$ in $L^1(F)$ where $F$ is a fundamental domain in $\tilde M $. Then for each $\gamma \in \pi_1(M)$, $\alpha_q(\gamma)$ is bounded in $\mathfrak g$ independent of $q$.
 \end{claim}
 \begin{proofconvvq}
 Choose coordinates 
 \[
 \varphi: [0,1] \times [0,1] \rightarrow \tilde M; \ \ (t,s) \mapsto \varphi(t,s)
 \]
 such that
 \begin{itemize}

 \item $\varphi(0,s) \in F$ and $\varphi(1,s)= \gamma \varphi(0,s)$ for every $s \in [0,1]$.
 \item For each $s \in [0,1]$, the curve $ t \mapsto \varphi(t,s)$ is a constant speed geodesic between $\varphi(0,s)$ and $\varphi(1,s)$.
 \item denote by  $\tilde v_q(t,s)$ the pullback under  $\varphi$ of $\tilde v_q$.
 \item $\tilde v_q(0,s)$ converges to $\tilde v(0,s)$ in $L^1([0,1])$.
 \end{itemize}
 
 All these assertions are easy to achieve, starting with last. Indeed, since $\tilde v_q \rightarrow \tilde v$ in $L^1$   we can choose, by Fubini's theorem, a path $\varphi(0,s)$ so that $\tilde v_q(\varphi(0,s))  \rightarrow \tilde v (\varphi(0,s))$ in $L^1([0,1])$. Then set $\varphi(1,s)= \gamma \varphi(0,s)$ and define $\varphi(t,s)$ as in the second bullet point by geodesic interpolation. The third assertion is just notation.
 
 Since $\tilde v_q \in W^{1,1}$,  we have $ t \mapsto \tilde v_q(t,s)$ is absolutely continuous for almost every $s \in [0,1]$ and
 \[
 \tilde v_q(\gamma \phi(0,s))-\tilde v_q(\phi(0,s)) =  \tilde v_q(1,s)-\tilde v_q(0,s)= \int_0^l\frac{\partial \tilde v_q(t,s)}{\partial t}dt.
 \]
Thus,
 \begin{eqnarray*}
 \int_0^l (Ad(\rho(\gamma)\tilde v_q(0,s)+\alpha_q(\gamma)-\tilde v_q(0,s))ds =
  \int_{[0,1] \times [0,1]}\frac{\partial \tilde v_q(t,s)}{\partial t}dtds.
 \end{eqnarray*}
 This implies that
 \begin{eqnarray*}
 \left |\int_0^l (Ad(\rho(\gamma)\tilde v_q(0,s)+\alpha_q(\gamma)-\tilde v_q(0,s))ds \right|
 \leq  \int_{[0,1] \times [0,1]} \left |d \tilde v_q(t,s)\right| dtds
 \end{eqnarray*}
 is uniformly bounded by a constant independent of $q$. Since $\tilde v_q(0,s)$ and thus also $Ad(\rho(\gamma) \tilde v_q(0,s)$ converge in $L^1$,  $\alpha_q(\gamma)$
 must be uniformly  bounded. 
 \end{proofconvvq}
 
 We now return to the proof of the theorem.  Using the boundedness of $a_q$ and following a diagonalization argument as in \cite[Theorems 4.3]{daskal-uhlen1}, we can arrange that $\alpha_q \rightarrow \alpha$ and $\tilde v_q \rightarrow \tilde v$  weakly in $BV_{loc}$ and strongly in $L^s_{loc}$ for  $s <2$.   This proves $(i)$ and $(ii).$ By construction $d \tilde v$ descends to $V$. Finally, we know $[V_q]=[dv_q]=[\alpha_q] \in H^1(M, ad(E_\rho))$, $[V_q]$ converges to $[V]$ and $\alpha_q $ converges to $\alpha $. It follows  $[V]= \alpha \in H^1(M, ad(E_\rho))$.
 \end{proof}
 
 \begin{theorem} \label{thm:limmeasures04rW} The analogue of Theorem~\ref{thm:limmeasures04r} holds for $W$. We write $dw=W.$
\end{theorem}

 \subsection{Remark on interpreting $v$, $w$ as sections of  affine bundles}\label{remarks} We end with explaining the term {\it affine representation} and interpreting $v$ as a section of an affine vector bundle. We continue with $G=SO^+(2,1)$, $\mathfrak g=so(2,1)$. Consider the semi-direct product  $Ad(G) \ltimes \mathfrak g$ consisting of pairs $(Ad(g), \alpha) \in Ad(G) \times \mathfrak g$ with group law
\[
(Ad(g_1), \alpha_1)\cdot  (Ad(g_2), \alpha_2)=(Ad(g_1g_2), Ad(g_1)\alpha_2+\alpha_1)
\]
and the homomorphism (called the linear part)
\[
L: Ad(G) \ltimes \mathfrak g \rightarrow Ad(G); \ \ L(Ad(g), \alpha)=Ad(g).
\]
The group $Ad(G) \ltimes \mathfrak g $ acts on $\mathfrak g$ by
\[
(Ad(g),\alpha) \in G \ltimes \mathfrak g \mapsto \{v \in \mathfrak g \mapsto Ad(g)v+\alpha\}.
\]
Given $\rho: \pi_1(M) \rightarrow G$,
then the action (\ref{affACt}) is nothing but a representation
\[
\hat \rho: \pi_1(M) \rightarrow Ad(G) \ltimes \mathfrak g
\]
 called an {\it infinitesimal  deformation of $Ad(\rho)$} (cf. \cite{gold} or \cite{dacing}).
 Such a representation $\hat \rho$ defines a flat {\it affine bundle} $\hat E=\tilde M \times_{\hat \rho}  E^{2,1}$ with linear part $ad(E_\rho)$ and an equivariant map $\tilde f$ as in (\ref{aqeqiov8}) defines a section $f: M \rightarrow \hat E$. Therefore, we have shown:
 \begin{theorem}\label{affsect} The primitives $\tilde v$, $\tilde w$ of $\tilde V$, $\tilde W$ define global sections $v$, $w$ of bounded variation of the affine vector bundles $\hat E_\rho=\tilde M \times_{\hat \rho}  E^{2,1}$,  $\hat E_\sigma=\tilde M \times_{\hat \sigma}  E^{2,1}$ such that $dv=V$, $dw=W$.
 \end{theorem}

\section{Lie algebra valued transverse measures}\label{prelimeart}
In  Section~\ref{measlieal} we define transverse measures on laminations with values in the Lie algebra of $SO(2,1)$.  An example is the Noether current on $M$ of an infinity harmonic map $u$. (cf. Section~\ref{sect:affbundlinfinite}). 
In Section \ref{sect:ruelle}
 we define the standard Lie algebra valued   transverse measure associated with a measured geodesic lamination. More precisely, if $\mu$  is a transverse measure on a geodesic lamination $\lambda$, we define the standard Lie algebra  valued transverse measure as $dw=B d\mu$ where $B$ denotes the generator of the geodesic flow of $\lambda$. 
In Section~\ref{exstthm} we prove the existence theorem that any Lie algebra valued  transverse measure  $W$ on a geodesic lamination induces a unique (real valued) transverse measure $\mu$ on the lamination.  Moreover $W=B d\mu$  (cf. Theorem~\ref{thmrepresentearth*} also Theorem~\ref{thmrepresentearth} in the introduction). In other words $W$ is the standard Lie algebra valued  transverse measure associated to $\mu$.
 Section~\ref{sect:length} contains the proof that the mass (total variation)  of  $W=B d\mu$ is equal to the $\mu$-length of the lamination  that supports $W$ (cf. Theorem~\ref{Length=totvar} in the introduction).

\subsection{Transverse measures with values in the Lie algebra}\label{measlieal} 
We will discuss certain types of measures with values in the Lie algebra $\mathfrak g=so(2,1)$.   
We fix a hyperbolic structure on $M$ given by a discrete, faithful representation $\sigma: \pi_1(M) \rightarrow G=SO^+(2,1)$. 
Let $\lambda$ denote a geodesic lamination on $M$  and $\tilde \lambda$ its lift to the universal cover $\tilde M$.  Recall that a  {\it plaque} is the closure in $\tilde M$ of a connected component of  $\tilde M \backslash \tilde \lambda$. Its image in $M$ will also be called a plaque. 

By a {\it flow-box}  for a geodesic lamination $\lambda$ (cf. \cite[Definition 7.2]{daskal-uhlen1}) we mean a bi-Lipschitz homeomorphism 
\begin{equation}\label{eqn:flowbox}
F: [a,b] \times [c,d]  \rightarrow  R \subset M; \ F=F(t,s)
\end{equation}
such that there exists a closed set $C \subset (c,d)$ of Hausdorff dimension 0 such that
\[
F^{-1}(\lambda)= [a,b] \times C.
\]
If we further assume that there exists a measure $\mu=\mu(s)$ on $[c,d]$ with support in $C$ and compatible with changes of flow-boxes, then $(\lambda, \mu)$ is a measured geodesic lamination.

In the following sections we will need the explicit construction of a flow-box (cf. \cite[Proposition 7.3]{daskal-uhlen1}), which we recall now:
 Let $f: [c,d] \rightarrow M$ be a Lipschitz transversal to the lamination. We can assume that $f$ is a geodesic and let
 \begin{equation}\label{defnB}
n: [c,d] \rightarrow TM
\end{equation}
 denote the Lipschitz function defined as follows. Let $C = \{ k \in [c,d]: f(k) \in \lambda \}$ and $n(k) \in T_{f(k)}(M)$ be  the unit tangent vector to the leaf of $\lambda$ through $f(k)$. By (cf. \cite[Lemma 7.1]{daskal-uhlen1}),  $n$ can be chosen to be  a Lipschitz function on $C$, which we can extend   to a Lipschitz function  on the interval by linearly interpolating on the plaques (which we still denote by $n$). Given any $k$, the continuity of $n$ implies that $n(s) \neq 0$ for $s$ close to $k$. For $s$ in an interval near $k$,
we set
\[
F(t,s) = exp_{f(s)}(t n(s)). 
\]
Note that
\begin{equation}\label{lipderF}
B(t,s):=\frac{\partial F(t,s)}{\partial t}=dexp_{f(s)}(t n(s))n(s)  
\end{equation}
is Lipschitz is $s$ and $C^\infty$ in $t$.

\begin{definition}\label{Def:LAMeasure} Let $\lambda$ be a geodesic lamination on $M$ (not necessarily measured). A {\it transverse measure on $\lambda$ with values in the Lie algebra} $\mathfrak g$ or simply a {\it Lie algebra valued transverse measure} is a function
\[
\tilde w: \tilde M \rightarrow  \mathfrak g^* \simeq \mathfrak g
\]
  (identify the Lie algebra and its dual via $(,)^\sharp$) satisfying the following properties:
 \begin{itemize}
 \item $(i)$ $\tilde w$ is constant on each plaque of $\lambda$
 \item $(ii)$ $\tilde w$ is of locally bounded variation
  \item $(iii)$ $\tilde w$ satisfies the equivariance relation
\begin{equation}\label{adimvf}
\tilde w(\gamma x)=Ad(\sigma(\gamma ))\tilde w(x)+ \alpha(\gamma), \ \ \forall \gamma \in \pi_1(M), x \in \tilde M
\end{equation}
where $\alpha$ is a 1-cocycle with values in $\mathfrak g$.
 \end{itemize}
 Under the above assumptions $d \tilde w$ descends to $M$ to define a distribution $dw$ of {\it finite mass} with values in $T^*M \otimes ad(E_\sigma)$ and support in $\lambda$. 
 We will further assume:
\begin{itemize}
 \item $(iv)$ $dw$ has values in $T^*M \otimes TM$ where $T^*M$ is considered as a subbundle of $ad(E_\sigma)$ via the map $\beta$ (cf. (\ref{defn:alphabetaf})). In the notation of Theorem~\ref{thm:limmeasures0}, 
$dw_{id}=dw.$
\end{itemize}
We also impose the positivity condition:
\begin{itemize}
 \item $(v)$   $*(\omega_{mc} \wedge dw)^{\sharp} $ is a (positive) Radon measure  on $M$. Here $ \omega_{mc}=dx \times x$ denotes the Mauer-Cartan form on $M$ (see comments after Theorem~\ref{thm:limmeasures0}). Note that the expression makes sense because $\omega_{mc}$ is smooth.
 \end{itemize} 
   {\it Often   we will denote the  Lie algebra valued transverse measure by $dw$ (instead of $\tilde w$).}
 \end{definition}
 \begin{remark} For people familiar with the topology literature, the above description of the measure $\tilde w$ is analogous to the one of transverse measures  via transverse cocycles. (See \cite[Section 7.2]{daskal-uhlen1} for the correspondence between transverse cocycles and functions of bounded variation.) But unlike the transverse cocycles  defined by Bonahon (cf. \cite{bonahon1}, \cite{bonahon3}),  our cocycles are not invariant under the action of $\pi_1(M)$: Instead they are twisted by the adjoint representation $Ad(\sigma)$. Also  recall that {\it we do not require that the support of} $dw$ {\it to be exactly} $\lambda$, {\it we only require that the support is  contained in} $\lambda$.
\end{remark}

%
%
 
 From the results of Section~\ref{sect:affbundlinfinite} we have shown that $W = dw$ associated to the Noether current of $u$ coming from the symmetries of the domain satisfies  properties $(i)$-$(v)$ of a Lie algebra valued transverse measure on the canonical lamination $\lambda$. 

\subsection{The standard Lie algebra valued transverse measure on a measured geodesic lamination}\label{sect:ruelle}

We start with a simple lemma which will be of importance later. 
Let $\{\gamma(t): t \in (-\epsilon, \epsilon)\}$  be a geodesic segment.
 We would like first to identify  the expression
$B_0:=\frac{d\gamma(t)}{dt} \times \gamma(t)$
which will appear in our formulas.
It follows from the geodesic equation for $\gamma$, that $B_0$ is constant along the geodesic.  
We can identify this constant geometrically:

\begin{lemma}\label{Lemma E.2} Suppose the geodesic segment $\gamma(t)$ is parametrized by arc length and is written as $e^{tB} X_0$ where $B \in \mathfrak g=\mathfrak s\mathfrak o(2,1).$  Then 
\[
\frac{d\gamma(t)}{dt} \times \gamma(t)= B.
\]
\end{lemma}
\begin{proof} Since $\gamma(t) = e^{Bt} X_0$,  
$\frac{d\gamma(t)}{dt}  = B e^{Bt} X_0$  
and 
\begin{eqnarray*}
B_0:=  \frac{d\gamma(t)}{dt}  \times \gamma(t) 
= e^{tB} B X_0 \times e^{tB} X_0
 = e^{tB} (B X_0 \times X_0) e^{-tB}
 \end{eqnarray*}
Since $B_0$ is  constant, 
\[
  0=  \frac{d}{dt} \left(e^{tB} (B X_0 \times X_0)e^{-tB} \right)= e^{tB} [B, BX_0 \times X_0]e^{-tB}.
\]
This means that $BX_0 \times X_0=bB$ is a multiple of $B$ and thus
\[
B_0 = e^{Bt} (bB) e^{-Bt}  = bB
\]
Since $\gamma$ is parametrized by arc length, $|B_0|=\sqrt 2.$
Since $\sqrt 2 \alpha_{X_0}$ is an isometry, 
\[
|\sqrt 2 \alpha_{X_0}(B)|=|B|=|\sqrt 2 \gamma'(0)|=\sqrt 2.
\]
Thus $b=1$.
\end{proof}

\begin{remark}The following is a reformulation of the lemma above: given a geodesic $\gamma(t)$ parametrized by arc length,  $B= \frac{d \gamma}{dt} \times \gamma(t)$ is the infinitesimal isometry  corresponding to the 1-parameter family of hyperbolic isometries with fixed axis $\gamma(t)$. Indeed,
\[
e^{sB}\gamma(t)=e^{sB}e^{tB}X_0=\gamma(t+s).
\]
This will be important when we discuss Thurston's earthquakes.
\end{remark}
 
  Let  $(\lambda, \mu)$ be a  measured geodesic lamination on $M$ and choose an orientation for $\lambda$ locally. For $X \in \lambda \cap R$, we denote by $B(X) \in \mathfrak g$  the element in the Lie algebra corresponding to the unit tangent field to the lamination consistent with the orientation. Note that $B$ is Lipschitz  along $\lambda$ and by the Lipschitz extension theorem it can be extended to a Lipschitz function in a neighborhood of the lamination \cite[Section 7.1]{daskal-uhlen1}. If $F=F(t,s): [a,b]\times [c,d] \rightarrow R \subset M$  is a flow-box of $\lambda$ as in (\ref{eqn:flowbox}),  we write 
 \begin{equation}\label{defB}
 \frac{\partial}{\partial t}=\frac{\partial F}{\partial t}, \ \  B=\frac{\partial} {\partial t} \times id.
 \end{equation} 
 $B$ is  {\it the generator of the geodesic flow of $\lambda$}.
  
  For a section $\phi \in \Omega^1(M, ad(E_\sigma))$ compactly supported in $R$, we define  
 \begin{equation*}\label{ruelle}
 T_\mu(\phi)=\int_c^d\int_{(a, b) \times \{s\}}(\phi, B)^\sharp d\mu(s). 
 \end{equation*}
   If we change the local orientation of the lamination, then both $B$ and $\partial /\partial t$  change sign and thus the value is unchanged. We can thus define $T_\mu(\phi)$ globally for any $\phi \in \Omega^1(M, ad(E_\sigma))$: Write $\phi=\sum_i \phi_i$ where $\phi_i$ is supported in a flow-box $R_i$ and set
   \begin{equation*}\label{ruelle2}
T_\mu(\phi)=\sum_i\int_{c_i}^{d_i}\int_{(a_i, b_i) \times \{s\}}(\phi_i, B)^\sharp d\mu(s). 
 \end{equation*}
This construction is similar to the scalar case (cf. \cite{sullivan} and \cite[Definition 8.1]{daskal-uhlen1}). However, the  difference with the scalar case is that here we do not need the extra assumption that $\lambda$ is oriented.  


\begin{theorem} \label{compatiblelam1} $T_\mu$ is a Lie algebra valued transverse measure. 
  \end{theorem}
  
  \begin{notation} 
  We will often  write $B d\mu$ instead of $T_\mu$ and call  it {\it the standard Lie algebra valued transverse measure} on $(\lambda, \mu)$.
 \end{notation}
    
  \begin{proof}
  The construction is  similar to \cite[Section 8]{daskal-uhlen1} so we will be brief. We first show that $T_\mu$ is closed. Assume that $\phi=d\xi$ where $\xi \in \Omega^0(ad(E_\sigma))$. Write $\xi= \sum_i\xi_i$ where $\xi_i$ is supported in the flow-box $R_i$. For $s \in [c_i, d_i]$ such that $(a_i, b_i) \times \{s\}$ is  a leaf of the lamination (otherwise $d\mu(s)=0$),
  \begin{eqnarray*}
  \lefteqn{\int_{(a_i, b_i) \times \{s\}}(d\xi_i, B)^\sharp =\int_{(a_i, b_i) \times \{s\}} (\frac{\partial \xi_i}{\partial t}, B)^\sharp dt}\\
  &=&\int_{(a_i, b_i) \times \{s\}} \frac{\partial }{\partial t}(\xi_i, B)^\sharp dt=0.
  \end{eqnarray*}
  In the second equation we used the geodesic equation along the lamination and in the last that $\xi_i=0$ on the boundary. It follows that $T_\mu(d\xi)=0$ and thus $T_\mu$ is closed.

  By the results of Sections~\ref{primV} and~\ref{remarks}, we can write $T_\mu=dw$ for  $w$ a section of an affine bundle with linear structure $ad(E_\sigma)$ (i.e condition $(iii)$ of Definition~\ref{Def:LAMeasure} holds). By construction, $T_\mu$ is supported on $\lambda$ (condition $(i)$ holds) and, since $\mu$ is a measure, $w$ is of finite variation (condition $(ii)$ holds). Because $B$ is tangent, condition $(iv)$ also holds. 
  
  We now show the positivity condition $(v)$. Fix a flow-box with coordinates $(t,s)$.  Write  $\omega_{mc}=(\frac{\partial}{\partial t} dt+ \frac{\partial}{\partial s} ds)\times id $. For $\chi \geq  0$ a test function,
  \begin{eqnarray}\label{explpos}
 \lefteqn{ (\omega_{mc} \wedge T_\mu)^\sharp (\chi) =  \sum_i\int_{c_i}^{d_i}\int_{(a_i, b_i) \times \{s\}}\chi_i (\omega_{mc}, B)^\sharp d\mu(s)} \nonumber\\
  &=& 2\sum_i\int_{c_i}^{d_i}\int_{(a_i, b_i) \times \{s\}}\chi_i  dtd\mu(s) \geq 0.
  \end{eqnarray}
 This completes the proof. 
  \end{proof}
     
  \begin{example} \label{exam:clgd}An important example of a Lie algebra valued transverse measure occurs when the lamination $(\gamma, b)=\{(\gamma_i, b_i)\}$ consists of a finite union of closed geodesics $\gamma_i$ and weights $b_i >0$. Assume $\gamma_i: [0, l_i] \rightarrow M$ is parametrized by arc length.
 First, we need to define a  special type of flow-box in a neighborhood $\mathcal O_i$ of each  geodesic called {\it parallel coordinates}:
 Assume $\gamma_i(t)=e^{tB_i}X_0$ and define 
  \[
  F_i: (0, 2\pi) \times (-\epsilon, \epsilon) \rightarrow \mathcal O_i=\mathcal O(\gamma_i) \subset M, \ \ F_i(t,s) = exp_{\gamma_i(t)}(s e_t)
  \]
  where  $e_t$ be the unit normal to $\gamma_i$ at $\gamma_i(t)$.
Note that $F_i(t,0)=\gamma_i(t)$, $s \mapsto  F_i (t,s)$ are geodesics parameterized by arc length perpendicular to $\gamma_i$ at $\gamma_i(t)$ and  $s \mapsto  F_i (t,s)$ are all equidistant nonintersecting curves. Furthermore, it is easy to see that the hyperbolic metric in the $(s,t)$ coordinate system is $ds^2 + \cosh s ^2 dt^2$. In particular, the metric is Euclidean on the geodesics $\gamma_i$. Set 
\[
T_b = \sum_i b_iB_i\delta_i(s)ds
\]
 where $\delta_i(s)$ denotes the $\delta$-function along $\gamma_i$.  
 $T_b$ is that standard Lie algebra valued transverse measure on $(\gamma, b)$. Note that, by (\ref{explpos}), condition $(v)$ of Definition~\ref{Def:LAMeasure} is equivalent to $b_i>0$.
   \end{example}

\subsection{Proof of Theorem~\ref{thmrepresentearth}}\label{exstthm} 
The purpose of this section is to provide a proof of the following theorem (stated also as Theorem~\ref{thmrepresentearth} of the introduction):

\begin{theorem}\label{thmrepresentearth*} Let $M$ be a closed hyperbolic surface and  $W=dw$  a Lie algebra valued  transverse measure with support in a  geodesic lamination $\lambda$ of $M$. Then  there exists a unique transverse measure $\mu$ with support in  $\lambda$  such that  $dw=Bd\mu$ where B is the generator of the geodesic flow of $\lambda$. Furthermore, 
$\mu=*(\omega_{mc} \wedge dw)^\sharp.$ Here, $\omega_{mc}=dX \times X$ denotes the Mauer-Cartan form on $M.$ (Note that multiplying a measure by a smooth form makes sense.)
  \end{theorem}  
In order to make the ideas of the proof more transparent we will consider first the case when the lamination consists only of closed geodesics and then deal with the general case.

\vspace*{.1cm}
    {\it {\bf Step 1:}  Theorem~\ref{thmrepresentearth*}  holds for  weighted closed  geodesics.}
    
    \vspace*{.1cm}

\noindent Assume the  closed geodesic $\gamma$ is given by $\gamma(t) = e^{Bt}X_0$ and $dw$ has support on $\gamma$. In parallel coordinates  centered at $\gamma$, choose $w(t,s)= 0$ for  $s <0$ and $w(t, s) = A$ for $s > 0$. In this case,  we have 
\begin{equation}\label{hlpedw}
dw = A\delta(s)ds 
\end{equation}
where $\delta(s)$ denotes the delta function along the geodesic. 
By property 
$(iv)$ of Definition~\ref{Def:LAMeasure}, $A_{\gamma(t)}=A_{e^{Bt}X_0}=A$ and by \cite[Proposition 3.5 $(v)$]{daskal-uhlen2},
\begin{eqnarray*}   
A_{\gamma(t)}=Ae^{Bt}X_0 \times e^{Bt} X_0= Ad e^{Bt} (( Ad e^{-Bt} (A) X_0 \times X_0)). 
\end{eqnarray*}
  There is a natural basis for the Lie algebra $\mathfrak g$ of $G=SO^+(2,1)$ compatible with the embedding of  $T_{X_0} \HH$ in $\mathfrak g$ via  cross product.  To make our computations, we select coordinates so that $X_0 =  (0, 0, 1)^\sharp$ and $B$ is the Lie algebra element
\[ 
              B = \{B_{ij}\} \ \mbox{where} \  B_{23}=B_{32} = 1 \ \ \mbox{and all other entries 0}.
\]
Let 
\[
B^\perp = {a_{ij}} \ \mbox{where} \ a_{13}  = a_{31}  = 1 \ \ \mbox{and all other entries 0.}
\]
Then $B$, $B^\perp$ span $T_{X_0} \HH$.
The normal direction  is
\[
\hat n = {d_{ij}}\ \mbox{where} \ d_{12}  = -d_{21}  = 1.
\]
First,
\begin{eqnarray*}   
 |A-A_{\gamma(t)}|=| A - Ad e^{Bt} ( Ad e^{-Bt} (A) X_0 \times X_0)| =
             |Ad e^{-Bt} (A) - Ad e^{-Bt} (A) X_0 \times X_0|.
\end{eqnarray*}
A straightforward computation with $A = bB + a B^\perp + z \hat n$
gives
\begin{equation}\label{matrcomp}
 |Ad e^{-Bt} (A) - Ad e^{-Bt} (A) X_0 \times X_0|=0 \  \ \forall t  \iff A=bB. 
 \end{equation}
 Indeed,\\  
$A=
\begin{pmatrix}
               & 0    & z      &a  \\                                
          &-z     &0            &b   \\                         
            & a       & b      &0                                
\end{pmatrix},
$ \ \ \ \ \ \ \ \ \ \ \ \ 
$
e^{tB}=
\begin{pmatrix}
               & 1    & 0      &0     \\                             
            &0    &\cosh t    &\sinh t    \\                     
            & 0       &\sinh t      &\cosh t                               
\end{pmatrix},
$\\
$
Ad(e^{-tB})A=
\begin{pmatrix}
               & 0    & z\cosh t -a\sinh t       &-z\sinh t +a\cosh t       \\                             
            &-z\cosh t +a\sinh t     &0     &b    \\                     
            & -z\sinh t  +a \cosh t      &b   &0                              
\end{pmatrix}
$, \\
$
Ad(e^{-tB})A-Ad(e^{-tB})AX_0 \times X_0=
\begin{pmatrix}
               & 0    & z\cosh t -a\sinh t       &0       \\                             
            &-z\cosh t +a\sinh t     &0     &0    \\                     
            & 0      &0   &0                              
\end{pmatrix}.
$
This implies formula (\ref{matrcomp}). Thus,
\[
dw = bB\delta(s)ds.
\]
Moreover, $b>0$ by property $(v)$ of Definition~\ref{Def:LAMeasure} (see also Example~\ref{exam:clgd}). 
Hence $dw$ is the standard Lie algebra valued transverse measure associated to $(\gamma, b)$. 
%
%
%
%
%

\vspace*{.1cm}
{\it {\bf Step 2:}  The general case.}
 \vspace*{.1cm}

The next proposition is the crucial ingredient in the proof of Theorem~\ref{thmrepresentearth*}. Before we state it, we record the following a simple observation: Consider a tangent vector $v \in T_X \HH$. Under the action of $G=SO^+(2,1)$ $gv \in T_{gX} \HH$.
Then
\begin{equation}\label{conjla}
gv \times gX=Ad(g)(v \times X).
\end{equation} 

\begin{proposition}\label{cruclemma}
Let $dw$ be  a Lie algebra valued  transverse measure on a geodesic lamination $\lambda$ and $F: [a,b] \times [c,d] \rightarrow R \subset M$ a flow-box as in (\ref{eqn:flowbox}). Then for any Lie algebra valued 1-form $\phi$ compactly supported  in $R$,
\begin{equation}\label{indepenoth}
dw(\phi)=dw((\phi, B)^\sharp B).
\end{equation}
Here $B=\frac{\partial F}{\partial t}$ is the generator of the geodesic flow of $\lambda$ (cf. (\ref{defB})).
\end{proposition} 

\begin{proof}

For each $(t,s) \in [a,b] \times [c,d]$  consider, as in Step 1, the  basis for the Lie algebra $\mathfrak g$ formed by the vectors 
$\{ B=B(t,s), B^\perp=B^\perp(t,s), \hat n= \hat n(t,s)\}$. Here $B(t,s)$ is  tangent to the leaf $s=cont$ of $\lambda$, $B^\perp(t,s)$ orthogonal to  $B(t,s)$  and $\hat n(t,s) $  normal to the span of $B$ and $B^\perp$. 

For fixed $\delta>0$ small, consider the map $g_\delta(X) = e^{B\delta}X$.  Equivalently, in the coordinates given by $F: [a,b]\times [c,d] \rightarrow M$ as in (\ref{eqn:flowbox}),  $g_\delta(t,s)=(t+\delta,s).$ Note that, because of (\ref{conjla}),
\begin{equation}\label{conjhatn}
\hat n(t+\delta,s)=Ad(e^{\delta B})\hat n(t,s).
\end{equation}
Since $g_\delta$ preserves the plaques,   
 \begin{equation}\label{invt1}
 dw(\phi) = dw({g_\delta}^*\phi)   
\end{equation}
for any Lie algebra valued 1-form $\phi$ supported in the flow-box. 
Set $\phi = \varphi \hat n$ for a 1-form $\varphi$.  Since $dw_{id}=dw$, (\ref{conjhatn}) and (\ref{invt1}) imply 
 \begin{equation}\label{invt2}
 dw((\varphi \circ g_\delta) Ad(e^{\delta B})\hat n)=dw({g_\delta}^*(\varphi \hat n))=dw(\varphi \hat n)=0.
  \end{equation}
  By the computation in Step 1, 
  \[
  Ad(e^{\delta B})\hat n=\sinh\delta B^\perp +\cosh\delta \hat n.
  \]
  Plug this into (\ref{invt2}) and use the fact that ${g_\delta}^*(\varphi)=\varphi \circ g_\delta$ is arbitrary, to obtain that for any compactly supported 1-form $\varphi$ in the interior of the flow-box,
 \begin{equation}\label{invt3}
  dw(\varphi B^\perp)=0.
 \end{equation}
Formulas (\ref{invt2}), (\ref{invt3}) imply (\ref{indepenoth}).
\end{proof}

It is now clear how to proceed with the proof of Theorem~\ref{thmrepresentearth*} in the general case. The  idea is the following:  By Proposition~\ref{cruclemma}  the only component of $dw$ is in the direction of $B$. Write $dw(t,s)=B(t,s)d\mu(t,s)$ for some measure $\mu(t,s)$. Since $w$ is constant on the plaques, $dw(t,s)$ is constant in $t$. Also,  $B(t,s)$ is constant in $t$ by the geodesic equation. Therefore $\mu=\mu(s)$ is a measure in the transverse direction satisfying $dw=Bd\mu$.  

Of course $dw$ is not a function but a distribution so some extra care is needed to make the rough argument above into a rigorous proof:
Consider $w=w(t,s)$ as a function defined on the rectangle $[a,b] \times [c,d]$. Since $w$ is constant on the plaques,  $w(t,s)=w(a,s)=:g(s)$ where $g$ is a  function of bounded variation on the interval $[c,d]$. Hence, there exists a Radon measure $\nu$  and a $\nu$-measurable Lie algebra valued function 
$U$ on $[c,d]$ with $|U|=1$ 
$\nu$-a.e, such that
\begin{eqnarray}\label{exprsull0}
\int_c^d (g(s), \frac{d\varphi(s)}{ds})^\sharp ds= \int_c^d (U(s), \varphi(s))^\sharp d\nu(s) 
\end{eqnarray}
for all Lie algebra valued smooth functions $\varphi$ compactly supported in the interior of $[c,d]$ (cf. Riesz representation theorem \cite[Chapter 6]{simon}). 

Let $\phi=\varphi_1dt+ \varphi_2ds$  compactly supported in $[a,b] \times[c,d].$ We claim
\begin{eqnarray}\label{exprsull}
dw(\phi)=\int_c^d\int_{(a, b) \times \{s\}} (U, \phi)^\sharp d\nu(s).
\end{eqnarray}
Indeed,
\begin{eqnarray*}\label{exprsull2}
dw(\phi)&=&\int_{[a,b] \times [c,d]}  \left(w(t,s), \left( \frac{\partial \varphi_1(t,s)}{\partial s}-\frac{\partial \varphi_2(t,s)}{\partial t}\right)\right)^\sharp dt \wedge ds\\
&=& \int_{[a,b] \times [c,d]}  \left(g(s),  \frac{\partial \varphi_1(t,s)}{\partial s} \right)^\sharp dt \wedge ds\\
&=& \int_c^d \left(g(s), \frac{d}{d s}\int_a^b \varphi_1(t,s) dt \right)^\sharp  ds\\
&=& \int_c^d  \left(U(s), \int_a^b \varphi_1(t,s)dt \right)^\sharp d\nu(s)\\
&=&\int_c^d\int_{(a, b) \times \{s\}} (U, \phi)^\sharp d\nu(s).
\end{eqnarray*}
In the second equality we used that $w$ is independent of $t$ and in the one before the last (\ref{exprsull0}). This proves (\ref{exprsull}).

Combining Proposition~\ref{cruclemma} and (\ref{exprsull0}),
\[
dw(\phi)=dw((\phi, B)^\sharp B)=\int_c^d\int_{(a, b) \times \{s\}} (\phi, B)^\sharp (U, B)^\sharp d\nu(s).
\]
By the geodesic equation $\frac{\partial}{\partial t}(U(s), B(s,t))^\sharp=0$. Set $d\mu(s) =(U(s), B(s,t))^\sharp d\nu(s)$ and this completes the proof of  Theorem~\ref{thmrepresentearth*}.

\subsection{Comment on convergence of laminations} The correspondence between measured geodesic laminations and Lie algebra valued transverse measures given in Theorem~\ref{thmrepresentearth} suggests the following identification between the topologies of the two spaces. First note that the space of Lie algebra valued transverse measures,
viewed as derivatives of functions of bounded variation, inherits a topology as measures. In the next theorem we show that this topology is the same as   the topology on the space of measured laminations  defined by Thurston (cf. \cite[Chapter 8]{thurston2}). For more details see \cite{aidan}.

\begin{theorem}\label{convlemma}
Let  ${(\lambda_k, \mu_k)}_{k=1,2,...}$,  $(\lambda, \mu)$ be  measured geodesic laminations and let $dw_k=B_kd\mu_k$, $dw=Bd\mu$ denote the corresponding Lie algebra valued  transverse measures on $\lambda_k$, $\lambda$ respectively.  Then $(\lambda_k, \mu_k) \rightarrow (\lambda, \mu)$ weakly in measure  if and only if $\tilde w_k \xrightharpoonup{BV_{loc}} \tilde w$. 
  \end{theorem} 
 \begin{proof} 
If $(\lambda_k, \mu_k) \rightarrow (\lambda, \mu)$ then their flow-boxes $F_k, F: [a,b] \times [c,d] \rightarrow M$ have the property that $F_k \rightarrow F$ and $\frac{\partial}{\partial t_k}=\frac{\partial{F_k}}{\partial t} \rightarrow \frac{\partial F}{\partial t} =\frac{\partial }{\partial t}$ in $C^0$ (cf. \cite[Corollary 1.12]{aidan}).
 The convergence properties of $F_k$ to $F$ imply that $B_k \rightarrow B$ uniformly.
Consider an open set $V \subset [a,b] \times [c,d]$ and a test function $\phi \in \Omega^1(M, ad(E_\rho))$ with compact support in the intersection of all $F_k(V)$. Then
\[
dw_k(\phi)=\int_V (B_k, \phi)^\sharp d\mu_k(s).
\]
 Since $F_k \rightarrow F$, $\frac{\partial}{\partial t_k} \rightarrow \frac{\partial }{\partial t}$ and $B_k \rightarrow B$ all converge uniformly   and  $\mu_k \rightarrow \mu$ weakly, we  have that $dw_k \rightarrow dw$ weakly, where $dw=Bd\mu$.
 
Conversely, suppose  that  $dw_k=B_kd\mu_k \rightharpoonup dw=Bd\mu$. Since $B_k$ and $B$ have norm $\sqrt 2$,   we have that the total variation of the measures $\mu_k$ is uniformly bounded away from infinity and zero. Thus (for example, see \cite[Theorem C]{aidan}), we can pass to a subsequence such that  $\mu_k \rightharpoonup \hat \mu$ for some transverse measure $\hat \mu$ on a geodesic lamination $\hat \lambda$. By the first part,  
$dw_k=B_kd\mu_k \rightharpoonup \hat B d\hat \mu$. 
We conclude that $\hat B=B$ and $\mu=\hat \mu$. 
 \end{proof}

 \subsection{Mass equals length}\label{sect:length}
We first recall the definition of the length of a measured lamination $(\lambda, \mu)$. Let $F_i: [a_i,b_i]\times [c_i, d_i] \rightarrow R_i \subset M$ be finitely many flow-boxes whose interiors cover  $\lambda$ and let $U=\cup int(R_i)$. Let $\psi_i$ be a partition of unity subordinate to the cover consisting on the interiors of $R_i$. We define 
\[
l_g(\lambda)=\sum_i \int_{c_i}^{d_i} \int_{[a_i, b_i] \times \{s\}} \psi_i \left|\frac{\partial F_i }{\partial t}\right |dt d\mu(s).
\] 
Also recall the definition of mass (or total variation) $|| \cdot||_{meas}$ from Definition~\ref{measrcurr}. The following  is Theorem~\ref{Length=totvar} of the introduction:
\begin{theorem}\label{Length=totvar*}Let $(\lambda, \mu)$ be a measured geodesic  lamination on the hyperbolic surface $(M,g)$ and $dw=Bd\mu$ a Lie algebra valued transverse measure. The mass (total variation) $||dw||_{meas}$ satisfies
\[
||dw||_{meas}=  2 l_g(\lambda).
\]
\end{theorem}
Note that in our construction before   we had normalized that $||dw||_{meas}=2$. 

\begin{proof}
We cover $\lambda$ by finitely many flow-boxes $F_i$ as before. 
 Let $ \frac{\partial}{\partial t_i}=\frac{\partial F_i}{\partial t}$ and
$\alpha_i=|\frac{\partial}{\partial t_i}|^{-1} (\frac{\partial}{\partial t_i})^\flat$ where $^\flat: TM \rightarrow T^*M$ denotes the musical isomorphism via the hyperbolic metric. Note that 
\begin{equation}\label{pull F_i}
|\alpha_i|=1 \ \ \mbox{and} \ \ F_i^*(\alpha_i)=\left|\frac{\partial F_i}{\partial t}\right |dt.
\end{equation}
We also choose a 1-form $\beta_i$ defined in the interior of $R_i$ so that $\beta_i\left(\frac{\partial}{\partial t_i}\right)=0$. This is equivalent to  
\begin{equation}\label{degbasiform}
(F_i^*\beta_i)\left(\frac{\partial}{\partial t}\right)=0
\end{equation}
i.e the restriction of $F_i^*\beta_i$ to the horizontal lines vanishes.
The forms $\alpha_i$ and $\beta_i$  form a local basis of $T^*M$. Furthermore, by definition, $\alpha_i$ and $\beta_i$ are orthogonal. We first prove the following claim:

\begin{claim}
Let $\phi \in \Omega^1(M, ad(E_\sigma))$, $\phi_{id}=\phi$  and write in the flow-box $R_i$
\begin{eqnarray*}
\phi=\xi_1 B+\xi_2B^\perp=(\xi_{11}^i\alpha_i+\xi_{12}^i\beta_i)B+(\xi_{21}^i\alpha_i+\xi_{22}^i\beta_i) B^\perp.
\end{eqnarray*}
Then $s(\phi) \geq \sqrt 2|\xi_{11}^i|$. 
\end{claim}
To see this claim, let $\phi=\xi \times id$ where $\xi_{id}=\xi \in \Omega^1(M, E)$. Then
\begin{eqnarray*}
\xi=\xi_1 e+\xi_2\hat e=(\xi_{11}^i\alpha_i+\xi_{12}^i\beta_i)e+(\xi_{21}^i\alpha_i+\xi_{22}^i\beta_i) \hat e
\end{eqnarray*}
where $\ e \times id=B$ and $ \hat e \times id=B^\perp$ and $e, \hat e$ is a local orthonormal basis of $TM$.
Since $s(\phi)=\sqrt 2s(\xi)$, it suffices to show
$s(\xi) \geq |\xi_{11}^i|$.

Write 
\begin{eqnarray*}
Q(\xi)^2=(\xi;\xi)
=|\xi_1^i|^2e \otimes e^\sharp+(\xi_1^i;\xi_2^i) (e \otimes {\hat e}^\sharp+\hat e \otimes e^\sharp)+|\xi_2^i|^2 \hat e \otimes{\hat e}^\sharp.
\end{eqnarray*}
Since $\{e, \hat e \}$ is an orthonormal basis, the largest eigenvalue of this matrix is 
\[
s(\xi^i)^2=\frac{1}{2}\left(|\xi_1^i|^2+|\xi_2^i|^2+\sqrt{(|\xi_1^i|^2-|\xi_2^i|^2)^2+4(\xi_1^i;\xi_2^i)^2}\right)\geq |\xi_1^i|^2\geq |\xi_{11}^i|^2.
\]
The right inequality follows because $\alpha_i$ and $\beta_i$ are orthogonal. 
The claim follows.

Let $\psi_i$ be a partition of unity subordinate to the cover consisting of the interiors of $R_i$. Let 
$\phi=(\xi_{11}^i\alpha_i+\xi_{12}^i\beta_i)B+(\xi_{21}^i\alpha_i+\xi_{22}^i\beta_i) B^\perp$ with compact support in $U$ and $s(\phi) \leq \sqrt 2$. 
Then 
\begin{eqnarray*}
dw(\phi) &=& 2\sum_i \int_{c_i}^{d_i} \int_{[a_i, b_i] \times \{s\}} F_i^*(\psi_i\phi, B)^\sharp d\mu(s)\\
&=& 2\sum_i \int_{c_i}^{d_i} \int_{[a_i, b_i] \times \{s\}} F_i^*(\psi_i)F_i^*(\xi_{11}^i\alpha_i+\xi_{12}^i\beta_i)d\mu(s)\\
&=&2\sum_i \int_{c_i}^{d_i} \int_{[a_i, b_i] \times \{s\}} F_i^*(\psi_i \xi_{11}^i) \left|\frac{\partial F_i}{\partial t}\right |dt d\mu(s)\\
&\leq&2\sum_i \int_{c_i}^{d_i} \int_{[a_i, b_i] \times \{s\}} F_i^*(\psi_i) \left|\frac{\partial F_i}{\partial t}\right |dt d\mu(s)\\\\
&=& 2l_g(\lambda).
\end{eqnarray*}
Here we used 
$(B,B)^\sharp=2$ in the second line, (\ref{pull F_i}), (\ref{degbasiform}) in the third  line and the claim above to deduce $|\xi_{11}^i| \leq 1$ in the forth line. 
Hence $||dw||_{meas} \leq 2l_g(\lambda).$ 

In order to see the reverse inequality
set   
$\phi= \sum_i \psi_i B \alpha_i$. Then $s(\phi)=\sqrt 2$ and  (\ref{pull F_i}) together with the fact that 
$(B,B)^\sharp=2$  imply 
\begin{eqnarray*}
 dw(\phi) &=& 2\sum_i \int_{c_i}^{d_i} \int_{a_i}^{b_i} F_i^*(\psi_i) \left|\frac{\partial F_i}{\partial t}\right |dt d\mu(s)\\
&=& 2l_g(\lambda).
  \end{eqnarray*}
Thus, $||dw||_{meas} \geq 2l_g(\lambda)$.
\end{proof}


\section{The transverse measure induced by best Lipschitz maps}\label{sect:KL}
Section~\ref{MainTHthmnot} contains the proof of Theorem~\ref{MMMainTHM} from the introduction (restated below as Theorem~\ref{mainthmnotdom*}).
In Section~\ref{sect:push}  we review the formula for the push-forward of a function of bounded variation by a Lipschitz map needed in the construction of the transverse measure on the target lamination.  
In Section~\ref{sect:targ} we discuss the measure $dv$ coming from the symmetries of the target. We use  the result in Section~\ref{sect:push}
to push $dv$ forward by the best Lipschitz map $u$ to obtain a  Lie algebra valued  transverse measure $dv^{\land}$ on $N$. We then apply the results of Section~\ref{prelimeart} to write  $dv^{\land}=B^{\land}d\mu^{\land}$ for a transverse measure $\mu^{\land}$ with support in the image $\lambda^{\land}=u(\lambda)$ of the canonical lamination. In Section~\ref{sect:relvw} we show that the measure $\mu^\land$ is the push-forward of the measure $\mu$ coming from $dw$ by the best Lipschitz map $u$ and provide a proof of Theorem~\ref{thm:1relvw} from the introduction (restated below as Theorem~\ref{thm:1relvw*}). 

\subsection{Transverse measures on the canonical lamination} \label{MainTHthmnot} We now prove the following theorem, stated also as Theorem~\ref{mainthmnotdom} in the introduction. \footnote {recall that we do not require that transverse measures are of full support. Here we differ  from the Thurston literature}

\begin{theorem}\label{mainthmnotdom*} Let $M$, $N$ be hyperbolic surfaces and let $u: M \rightarrow N$ be an infinity harmonic map in a given homotopy class. Let $\lambda$ be the canonical lamination corresponding to the homotopy class and  $W=dw$  the measure on $M$ coming from the Noether current of $u$.  Then  
 $dw$ is a Lie algebra valued transverse measure on $\lambda$ and $dw=B d\mu$, where $\mu=2|S| :=2\lim_{p \rightarrow \infty} Tr(Q(\kappa_p du_p)^p)$.  
 \end{theorem}


\begin{proof}Theorems~\ref{thm:limmeasures0}  and~\ref{thm:limmeasures04rW} immediately imply that  $dw$ is a Lie algebra valued transverse measure on $\lambda$.   Theorem~\ref{thmrepresentearth*} implies that $dw=B \mu ds$ for a transverse measure $\mu$ such that $\mu= *(\omega_{mc}  \wedge dw)^\sharp$. Theorem~\ref{thm:limmeasures0}$(iv)$  implies that the latter is equal to $2|S|$.
 \end{proof}
 
\begin{corollary}\label{MMMainTHM*}Let $M$, $N$, $u$ and $\lambda$ as before. Assume that such that $u=\lim_{p \rightarrow \infty} u_p$  for a sequence  $u_p$ of $J_p$-minimizers in the same homotopy class. Then after passing to a subsequence  and after normalizing by positive constants $\kappa_p$, $|S_{p-1}| := Tr(Q(\kappa_p du_p)^p)$ converges weakly to a transverse measure 
$|S|$ on the canonical lamination $\lambda$ associated to the hyperbolic metrics on $M$, $N$ and the homotopy class of $u$.
\end{corollary}
\begin{proof}
Since $\mu=*(\omega_{mc} \wedge dw)^\sharp =2|S|$ is a transverse measure on $\lambda$, the proof follows immediately from Theorem~\ref{mainthmnotdom*}.
\end{proof}
 
\subsection{The push-forward of a BV function}\label{sect:push} So far we have discussed the Lie algebra valued measure $W=dw$ on $M$ induced from the Noether current of $u$ on the domain. Before we discuss the Noether current $V=dv$ on the target, we review the push-forward of a BV function by a Lipschitz map.  In order to give some motivation for the construction, first recall the following simple consequence of the change of variables formula. 

Let $g: \Omega \rightarrow \Omega^{\land}$ be a diffeomorphism between two open domains in $\R^n$ and 
$v:  \Omega \rightarrow \R^N$ an $L^1$ function. Define $v^\land=g_*(v)=v \circ g^{-1}: \Omega^{\land} \rightarrow \R^N.$ Let   $\psi \in C^\infty(\Omega^{\land})$ and $\phi \in \Omega^1(\Omega^{\land})$  both with support in $\Omega^{\land}.$  Then the change of variables formula implies
\[
 \int_{\Omega^{\land}} v^{\land}*\psi =\int_{\Omega} v g ^*(*\psi)\ \mbox {and} \ \int_{\Omega^{\land}} v^{\land}d\phi=\int_\Omega  v d(g^*\phi).
\]
 We need to apply the above formula for $g=u$ a best Lipschitz map which is not necessarily a difeomorphism. 
 
Let $g: \Omega \rightarrow \Omega^{\land}$ be a Lipschitz, proper  map between two domains in $\R^n$ of Lipschitz constant $L$. Let $v:  \Omega \rightarrow \R^N$  be a  function of bounded variation.   
 Define the push-forward of  $v$ by $g$
\begin{equation}\label{eqndefnpushforw}
 v^{\land}(y):=  g_*(v)(y)= \sum_{x \in  g^{-1}(y)}deg_{x}(g){v}(x).
 \end{equation}
 We set $v^{\land}=0$ outside the image of $g$. Here $deg_x(g)$ denotes the topological degree of $g$. Note that this function is defined outside the set $ g(E)$ where $E$ is the set of points where $ g$ is not differentiable or a critical point of $ g$. This set is of measure zero by \cite[Lemma 2.73]{ambrosio}.   The next lemma characterizes uniquely $v^{\land}$ and $dv^{\land}$.
 
 \begin{lemma}\label{pushforamb} Let $v$, $g$ and $v^\land$ as above. Then  $v^{\land}: \Omega^{\land} \rightarrow \R^n$ is BV and 
 \begin{equation}\label{normpushforw}
 ||dv^{\land}|| \leq C ||dv||
 \end{equation}
 where $C$ depends only on the Lipschitz constant of $g$ and the geometry.  
 Moreover, for  $\psi$   a smooth $\R^N$-valued  function and $\phi$  a smooth $\R^N$-valued  1-form both with support in $\Omega^{\land}$
 \begin{equation}\label{areaform}
 \int_{\Omega^{\land}} (v^{\land},*\psi)^\sharp=\int_{\Omega} (v, g^*(*\psi))^\sharp\ \mbox{and} \ \int_{\Omega^{\land}} (v^{\land}, d\phi)^\sharp=\int_\Omega  (v, dg^*(\phi))^\sharp.
 \end{equation}.
 \end{lemma}
 \begin{proof}
 The fact that $v^{\land}$ is of bounded variation and formula (\ref{normpushforw}) are contained in \cite[Theorem 3.16]{ambrosio}. The first formula in (\ref{areaform}) follows from the area formula (cf. \cite[(3.13)]{ambrosio}). The second follows from the first as follows: We write
 $d\phi=*\psi$ and apply the first formula
 \begin{eqnarray*}
 \int_{\Omega^{\land}} (v^{\land}, d\phi)^\sharp=\int_{\Omega^{\land}} (v^{\land}, *\psi)^\sharp=\int_{\Omega} (v, g^*(*\psi))^\sharp=\int_\Omega  (v, dg^*(\phi))^\sharp.
 \end{eqnarray*}
\end{proof}

\begin{remark}
Note that in the above  $||.||=||.||_{BV}$ is used to denote the usual total variation. If we use instead $||.||_{meas}$,  $C$ gets replaced by  another constant. This has no effect on our results.
\end{remark}

 \subsection{The Lie algebra valued  transverse measure on the target}\label{sect:targ} 
 Let  $\tilde u: \tilde M=\HH \rightarrow \HH=\tilde N$ denote the lift of $u$ to the universal cover intertwining the representations $\sigma$ and $\rho$ defining the hyperbolic structure on $M$ and $N$. Let $\tilde v: \tilde M=\HH \rightarrow \mathfrak g$ be, as in Section~\ref{primV}, the primitive of the lift of $V$.   By the compactness of the surfaces $M$ and $N$,  $\tilde u$ is proper. Set 
 \begin{equation}\label{def:hatv}
\tilde v^{\land}(y):= \tilde u_*(\tilde v)(y) 
\end{equation}
  the push-forward of $\tilde v$ as in (\ref{eqndefnpushforw}). 
    
  \begin{proposition}\label{propi-iv}$\tilde v^{\land}$ satisfies properties $(i)$-$(v)$ of Definition~\ref{Def:LAMeasure}.
 \end{proposition}
 \begin{proof}
 We first check property $(i)$, i.e $\tilde w$ is constant on each plaque of $\lambda$.  Fix fundamental domains $\Omega$ and $\Omega^{\land}$ in $\tilde M$, $\tilde N$ and denote by $v$, $v^{\land}$ the restriction of  $\tilde v$, $\tilde v^{\land}$.
  We first show that $v^{\land}$ is constant on the plaques of the lift $\tilde \lambda^{\land}$ of $\lambda^{\land}$ to $\HH$. Let $\phi$ be a smooth 1-form, as before, be supported in an open subset of $\Omega^{\land} \backslash  \tilde \lambda^{\land}$. Since $\tilde u$ sends $\tilde \lambda$ to $\tilde \lambda^{\land}$, $u^*(\phi)$ is supported in $\tilde M \backslash  \tilde \lambda$.  Since $\int_M (v, du^*(\phi))=0$, (\ref{areaform}) implies
 $dv^{\land}=0$ distributionaly. This proves property $(i)$.   
 Property $(ii)$ stating that $\tilde w$ is of locally bounded variation  follows from Lemma~\ref{pushforamb}.
Lemma~\ref{propiii}  below implies that  $\tilde w$ satisfies the equivariance relation, hence  property $(iii)$ holds.

We now prove property $(iv)$, i.e $dw_{id}=dw.$ Let 
  $\phi \in \Omega^1(N, ad(E_\rho))$ compactly supported in  $\Omega^{\land}$.  By~(\ref{areaform}),
  \[
  \int_{\Omega^{\land}} (v^{\land}, d(\phi_{id}))^\sharp=\int_\Omega(v, du^*(\phi_{id}))^\sharp=\int_\Omega(v, d(u^* \phi)_u)^\sharp.
\]
Here we used $u^*(\phi_{id})=(u^* \phi)_u$.
By Theorem~\ref{thm:limmeasures0}($v$) $dv_u=dv$, hence
\[
  dv^{\land}(\phi_{id})=dv((u^* \phi)_u)=dv(u^* \phi)=dv^{\land}(\phi)
\]
completing the proof of $(iv)$.

We finally check $(v)$, i.e   $*(\omega_{mc} \wedge dw)^{\sharp} $ is a (positive) Radon measure  on $M$.       Pick a sequence  $u_p \rightarrow u$ uniformly and $v_q \rightarrow  v$  weakly in $BV_{loc}$ as $p \rightarrow \infty$ ($q \rightarrow 1$).  Call $v_q^\land=u_{p*}(v_q)$. First notice that, since $v_q \rightarrow v$ in $L^1$ and thus almost everywhere, (\ref{eqndefnpushforw}) implies that $v^\land_q \rightarrow v^\land$ almost everywhere as well.  Next note that, by (\ref{normpushforw}), $||dv_q^\land||$ are uniformly bounded. Thus, after passing to a further subsequence, $v^\land_q $ converge locally weakly in BV. Since $v_q^\land \rightarrow v^\land $ almost everywhere, we conclude that 
\begin{equation}\label{convvland}
v_q^\land \rightarrow v^\land \ \ \mbox{locally weakly in BV}.
\end{equation}
Next pick $\eta \geq 0$ a smooth function in $N$. Then
\[
<dv_q; (\eta \circ u_p) du_p \times u_p>=2\int_M (\eta \circ u_p) TrQ^p(du) \geq 0.
\]
Thus, 
\begin{eqnarray*}
0 &\leq& <dv_q; (\eta \circ u_p) du_p \times u_p>= -\int_M (v_q, d*({u_p}^*(\eta \omega_{mc})))^\sharp\\
&=&-\int_N (v^\land_q, d*\eta \omega_{mc})^\sharp =<dv^\land_q; \eta \omega_{mc}>.
\end{eqnarray*}
 In the first equality we used ${u_p}^*(\omega_{mc})=du_p \times u_p$ and in the second  (\ref{areaform}). 
By taking limits and using (\ref{convvland}), we obtain
\[
<dv^\land; \eta \omega_{mc}> \geq 0.
\]
Since this is true for any smooth function $\eta \geq 0$ and, on the support of $dv$, the map $u$ is a homeomorphism, property $(v)$ follows.
 \end{proof}
 
 \begin{lemma}\label{propiii} Identify $\tilde M$ and $\tilde N$ with $\HH$ and the action of $\pi_1(M) \simeq \pi_1(N)$ with the corresponding images in $SO^+(2,1)=PSL(2, \R)$. Then 
 \[
  \tilde v^{\land}(\rho(\gamma)y)=Ad(\rho(\gamma))\tilde v^{\land}(y)+\alpha(\gamma)
 \]
 satisfies the same equivariance as $\tilde v$.
 Thus, $d\tilde v^{\land}$ descends to a  distribution $dv^{\land}$ on $N$ with values on $T^*N \otimes ad(E_\rho)$.
 \end{lemma}
 \begin{proof} 
 \begin{eqnarray*}
 \tilde v^{\land}(\rho(\gamma)y)&=&\sum_{x \in \tilde u^*(\rho(\gamma)y)}deg_{\tilde u}(x){\tilde v}(x) \nonumber\\
 &=&\sum_{x' \in \tilde u^*(y)}deg_{\tilde u}(\sigma(\gamma)x'){\tilde v}(\sigma(\gamma)x')\\
 &=&Ad(\rho(\gamma))\sum_{x' \in \tilde u^*(y)}deg_{\tilde u}(\sigma(\gamma)x'){\tilde v}(x')+\alpha(\gamma)\sum_{x' \in \tilde u^*(y)}deg_{\tilde u}(\sigma(\gamma)x')\\
 &=&Ad(\rho(\gamma))\sum_{x' \in \tilde u^*(y)}deg_{\tilde u}(x'){\tilde v}(x')+\alpha(\gamma)\sum_{x \in \tilde u^*(\rho(\gamma)y)}deg_{\tilde u}(x)\\
 &=&Ad(\rho(\gamma))\tilde v^{\land}(y)+\alpha(\gamma).
 \end{eqnarray*} 
 In the second to last equality we used the fact that  Fuchsian groups consist of orientation preserving isometries and the last equality holds because $deg(u)=1.$ 
 \end{proof}
 
 \begin{corollary}\label{asse1-3} Let $\lambda^{\land}=u(\lambda)$. Then:
\begin{itemize}
\item $(i)$ $dv^{\land}$ is a  Lie algebra valued  transverse measure on $N$ with support in $\lambda^{\land}$ (more specifically with support in the image of the support of $dv$ under $u$) 
\item $(ii)$ There exists a transverse measure $\mu^\land$ on $\lambda^{\land}$ such that $dv^{\land}=B^{\land}d\mu^{\land}$, where $B^{\land}$ is the generator of the geodesic flow of $\lambda^{\land}.$ Moreover, $\mu^{\land}=*(\omega_{mc}^{\land} \wedge dv^{\land})^\sharp$ where $\omega_{mc}^{\land}$ is the Mauer-Cartan form on $N$.
\end{itemize} 
 \end{corollary}
 
 \begin{proof}
 Combine Proposition~\ref{propi-iv} with Theorem~\ref{thmrepresentearth}. 
 \end{proof}

  \subsection{The relation between $dv$ and $dw$} \label{sect:relvw}

\begin{lemma}\label{wconvneg} Assume that $u_p \rightarrow u$ is a $p$-approximation by minimizers of $J_p$. Then $du_p \times u_p \rightharpoonup du \times u$ and $du_p \times du_p \rightharpoonup du \times du$ weakly in $L^s$ for as $s \geq 1$.
 \end{lemma} 
 \begin{proof} By Theorem~\ref{thm:existence}, $u_p \rightarrow u$ strongly in $C^0$ and weakly in $W^{1,s}$ for all $s$. Thus, $du_p \xrightharpoonup{L^s} du$ and  $du_p \times u_p \xrightharpoonup{L^s} du \times u$. This also implies that $du_p \times du_p=d(du_p \times u_p) \xrightharpoonup{L^s_{-1}} d(du \times u)=du \times du$  for all $s$.
 On the other hand, since $du_p$ (and thus also $du_p \times du_p$)  is uniformly bounded in $L^s$ for all $s$, after passing to a subsequence, $du_p \times du_p \xrightharpoonup{L^s} \Psi$, where $\Psi$ is a section of $\Lambda^2(T^*M) \otimes ad(E_\rho)$. It follows $\Psi=du \times du$ and completes the proof. 
 \end{proof}
 
 Before we proceed with the next proposition we make the following observation:
 Let $\Phi$ be a 1-form with values in  $ad(E_\rho)$ such that both $\Phi$ and $d\Phi$ are in $L^\infty$. We define the distribution $*(\Phi \wedge dv)^\sharp$ acting on a smooth  function $\chi$ on $M$ by
 \begin{equation}\label{anzfor}
 *(\Phi \wedge dv)^\sharp (\chi)=\int_M (d\Phi \wedge v \chi  )^\sharp -\int_M (\Phi v \wedge d\chi )^\sharp.
 \end{equation}
 
 By approximating with smooth sections it is easy to check that $*(\Phi \wedge dv)^\sharp$ is a well defined as a distribution (cf. \cite{anzellotti}).

 \begin{lemma}\label{comprdvdw} Assume that $u_p \rightarrow u$ is a $p$-approximation by minimizers of $J_p$ and $v_q$ are the dual sections such that $dv_q=*Q(du_p)^{p-2} du_p \times u_p \rightharpoonup dv$ weakly as measures. Then  
 \[
 *(du_p \times u_p  \wedge dv_q)^\sharp \rightharpoonup *(du \times u \wedge dv)^\sharp
 \]
 weakly as measures. 
  \end{lemma}
 \begin{proof} We are going to apply (\ref{anzfor}) for $\Phi=du \times u$ where $u$ is an infinity harmonic map.  In this case 
 $d\Phi=du \times du$ so both $\Phi$ and $d\Phi$ are in $L^\infty$. Thus $*(\Phi \wedge dv)^\sharp$ is well defined. 
 
 Let $\chi$ a smooth  function on $M$. Set 
 $\Phi_q=du_p \times u_p$. Lemma~\ref{wconvneg} implies that $\Phi_q \xrightharpoonup{L^s} \Phi $ and $d\Phi_q \xrightharpoonup{L^s} d\Phi $ for any $s>2$. 
 The assumption that $v_q $ converge weakly in $BV$ implies that $v_q $ converge strongly in $L^r$ for any $r<2$. Let $s >2$ be the conjugate to $r$, i.e. $1/s+1/r=1$. Then $\int_M (d\Phi_q \wedge v_q \chi  )^\sharp \rightarrow \int_M (d\Phi \wedge v \chi  )^\sharp$ and similarly for the term involving $\Phi$. This completes the proof of the lemma.
  \end{proof}
 
\begin{corollary}\label{cor:psh} $u_*(\mu)=\mu^\land$, i.e $\mu^\land$ is the push-forward of $\mu$ by $u$.
 \end{corollary}
  \begin{proof}
  First note that 
  \begin{eqnarray*}
 *(du_p \times u_p  \wedge dv_q)^\sharp= 2(du_p  \wedge *S_{p-1})^\sharp=2(du_p; S_{p-1})^\sharp=2|S_{p-1}|.
 \end{eqnarray*}
 and thus, by taking limits,
  \begin{eqnarray}\label{ecl1}
 *(du \times u  \wedge dv)^\sharp= 2|S|=\mu.
 \end{eqnarray}
  Let $\chi^\land \in C^\infty(N)$ and  set $\chi=u^\star(\chi^\land)=\chi^\land \circ u \in C^\infty(M)$. Then
  \begin{eqnarray}\label{ecl2}
 \lefteqn{*(du \times u  \wedge dv)^\sharp(\chi)} \nonumber\\
 &=& \int_M (du \times du \wedge  u^*(\chi^\land) v)^\sharp -\int_M (du \times u  \wedge du^*(\chi^\land) v )^\sharp   \nonumber\\
 &=& \int_M (u^*(d\omega_{mc}^{\land} \ \chi^\land) \wedge  v)^\sharp -\int_M (u^*(\omega_{mc}^{\land} \wedge d\chi^\land)v )^\sharp  \\
  &=&  \int_N (d\omega_{mc}^{\land} \wedge  v^\land \chi^\land )^\sharp -\int_N (\omega_{mc}^{\land} \wedge v^\land d\chi^\land  )^\sharp \nonumber\\
 &=& *(\omega_{mc}^{\land}   \wedge dv^\land )^\sharp (\chi^\land) = \mu^\land(\chi^\land) \nonumber.
 \end{eqnarray}
 In the above, the first and the forth equalities are applications of (\ref{anzfor}). 
In the third, we use (\ref{areaform}) and 
 in the last equality  Corollary~\ref{asse1-3}$(ii)$.
 Combining (\ref{ecl1}) with (\ref{ecl2}), we obtain
  \begin{eqnarray*}
 \mu(u^*(\chi^\land))=\mu(\chi)=2|S|(\chi)=\mu^\land(\chi^\land).
 \end{eqnarray*}
 Thus,  $u_*\mu=\mu^\land$.
 \end{proof} 
    
    We are now ready to prove the following theorem (Theorem~\ref{thm:1relvw} from the introduction).
 \begin{theorem}\label{thm:1relvw*} Let $M$,$N$, $u$, $\lambda$  as in Theorem~\ref{mainthmnotdom*}. Then
  \begin{equation*}
dv=B^\land(u) d\mu
\end{equation*}
where $B^\land$ is the generator of the geodesic flow of the geodesic lamination $\lambda^{\land}=u(\lambda)$ and $\mu=2|S|=*(\omega_{mc}, dw)^\sharp$. 
 \end{theorem}
  Let $F: [a,b] \times [c,d] \rightarrow R$, $F^\land: [a^\land,b^\land] \times [c^\land,d^\land] \rightarrow R^\land$ be flow-boxes for $\lambda$, $\lambda^{\land}$ in $M$, $N$ with $u(R) \subset R^\land$. Write  $(t^\land, s^\land)=u(t, s)$ and let $\phi^\land$  be a Lie algebra valued 1-form with support in  $R^\land$. Then 
 \begin{eqnarray*}
 dv(u^*\phi^\land)&=& -\int_M(v, d(u^*\phi^\land))^\sharp
 = -\int_N(v^\land, d\phi^\land)^\sharp\\
 &=& -\int_{c^\land}^{d^\land}\int_{a^\land}^{b^\land}(B^\land(t^\land,s^\land), d\phi^\land (t^\land, s^\land))^\sharp d \mu^\land(s^\land)\\
 &=& -\int_c^d\int_a^b(B^\land(u(t, s)), d(u^*\phi^\land)(t, s))^\sharp d\mu(s).
 \end{eqnarray*}
In the second equality we used  (\ref{areaform}), in the third  Corollary~\ref{asse1-3} and in the last equality we used change of variables together with Corollary~\ref{cor:psh}. Since the above equality is true for all $\phi^\land$ and $u$ is a homeomorphism on the support of $dv$, $dv=B^\land(u) d\mu$. This completes the proof.

We end the section with proving a precise formula relating  the Noether currents $dv$, $dw$ and the best Lipschitz map $u_p$. First recall that by Proposition~\ref{prop:relvw} the {\it renormalized tensors} $dv_q$ and $dw_q$ are related by
\[
-2T_q=(*dv_q \otimes U_p \times u_p)^\sharp, \ \ \ *dw_q=-2T_q \times id.
\]
Here $U_p=\kappa_p du_p \rightarrow L^{-1} du \times u$ where $L$ is the best Lipschitz constant.
$T_q$ denotes the energy-momentum tensor for the $p$-approximation $u_p$ of the best Lipschitz map. 

Because of the discontinuity of $du$, this formula relating $dv_q$, $dw_q$ and $u_p$  is difficult to pass to the limit. However, we are able to prove the  limiting formula useing  Theorems~\ref{mainthmnotdom*} and~\ref{thm:1relvw*}:

\begin{corollary}\label{ncdomtargu*} The Noether currents $dw$, $dv$ and the best Lipschitz map $u$ are related by
\[
-2T=L^{-1}(*dv \otimes du \times u)^\sharp, \ \ \ *dw=-2T \times id.
\]
Here $T$ is the limiting energy-momentum tensor and $L$  the best Lipschitz constant.
\end{corollary}
 \begin{proof} First note that, since $du$ is continuous on the support of $dv$, the expression above is well defined. More precisely, on the support of $\lambda$, write $du=L \partial_t \otimes  {\partial_t}^\land(u)$ where $\partial_t$ and $\hat \partial_t$ are the unit tangent vectors to the leaves of the laminations $\lambda$ and $u(\lambda)$.
 By definition
 \[
 B=\partial_t \times id \ \ \mbox{and} \ \ B^\land(u)={\partial_t}^\land(u) \times u.
 \]
 Using Theorem~\ref{thm:1relvw*},
 \begin{eqnarray*}
L^{-1}(*dv \otimes du \times u)^\sharp 
 &=& L^{-1}(*B^\land(u) d\mu \otimes du \times u)^\sharp \\
 &=& (*{\partial_t}^\land(u) \times u \ d\mu \otimes (\partial_t \otimes  {\partial_t}^\land(u))\times u)^\sharp\\
 &=& 2(*{\partial_t}^\land(u)  \ d\mu \otimes (\partial_t \otimes  {\partial_t}^\land(u)))^\sharp\\
 &=& 2*\partial_td\mu.
 \end{eqnarray*}
 Since $dw=*T \times id$,
 \[
 L^{-1}(*dv \otimes du \times u)^\sharp \times id=2*dw=-2T \times id.
 \]
 This completes the proof of the corollary.
\end{proof}

 \section{Lie algebra valued transverse measures are infinitesimal earthquakes}\label{lamearth}
Consider the  Teichm\"uller space  of a closed surface $M$ of genus at least two viewed as the quotient space of the space of hyperbolic metrics  by the group of diffeomorphisms homotopic to the identity.  The  Teichm\"uller space can also be identified with the space $ \mathcal R (M)$ of conjugacy classes of discrete, faithful representations of $\pi_1(M) \rightarrow PSL(2,\R)=SO^+(2,1)$. For $\sigma \in \mathcal R (M)$ the tangent space $T_\sigma \mathcal R (M)$  is equal to $H^1(\pi_1(M), \mathfrak g_{Ad(\sigma)})$ (cf. \cite{gold}). By the De Rham Theorem~\ref{cohoisp}, $T_\sigma \mathcal R (M)$ can also be identified with $ H^1(M, ad(E_\sigma))$ of closed $ad$-valued 1-forms modulo exact ones.  Furthermore,  there is a symplectic 2-form $\Omega$ on $\mathcal R(M)$   given by 
 \[
 \Omega([\phi], [\psi])= \frac{1}{2}\int_M \phi \wedge \psi.
 \]
 Here $\phi, \psi$ are closed 1-forms with values in the flat Lie algebra bundle $ad(E_\sigma)$ and $ [\phi], [\psi]$ denotes their cohomology classes in 
 $H^1(M, ad(E_\sigma)) \simeq T_\sigma  \mathcal R (M).$ The form $\Omega$ is equal (up  to a constant multiple) to the Weil-Petersson K\"ahler form on Teichm\"uller space  (cf. \cite{ahlfors} and \cite{gold}).
 
 Given a hyperbolic structure  $\sigma$ on $M$, a  simple closed geodesic $\gamma$   and $t \in \R$, we can deform  $\sigma$ by cutting the surface along $\gamma$, turn left-hand side of $\gamma$ in the positive direction at distance $t$ and reglueing back. We  denote the new surface by $E^t_\gamma(\sigma)$. As $t$ varies, we obtain a continuous path in  Teichm\"uller space called the Fenchel-Nielsen deformation of $\sigma$ along $\gamma$.  This construction extends to include all measured laminations $(\lambda, \mu)$ by defining 
 \[ 
 E^t_{(\lambda, \mu)}(\sigma)=\lim_i E^t_{(\gamma_i, \nu_i)}(\sigma)
 \]
 where $(\gamma_i, \nu_i) \rightarrow (\lambda, \mu)$ is a sequence of finite laminations converging in the space of projective measured laminations. This is called the {\it earthquake flow} in Teichm\"uller space and it is shown to depend smoothly on $\sigma$ and  $t$ (cf. \cite{ker1}, \cite{ker2}).
 
 The tangent vector to Teichm\"uller space at $\sigma$ given by
 $
 \frac{d(E^t_{(\lambda, \mu)}(\sigma))}{dt}\Big|_{t=0}
 $
 is called the {\it infinitesimal earthquake} at $\sigma$ (with respect to $(\lambda, \mu)$). An important property is that the infinitesimal earthquake is dual with respect to the   symplectic form $\Omega$ to the derivative of the length function of $(\lambda, \mu)$ at $\sigma$. In other words, consider
 \[ 
\mbox{length}_{\lambda, \mu} :  \mathcal R(M) \rightarrow \R
\] 
 as a function on Teichm\"uller space associating to $\sigma$ the length of $(\lambda, \mu)$ with respect to $\sigma$. Then for any closed 1-form $\phi$ with values in $ad(E_\sigma)$
\begin{eqnarray}\label{symplgradl}
(d\mbox{length}_{\lambda, \mu})_\sigma([\phi])=\Omega\left( dE^t_{(\lambda, \mu)}(\sigma)/dt \big|_{t=0},[\phi] \right).
\end{eqnarray}
 Here the cohomology class $[\phi] \in H^1(\pi_1(M), {\mathfrak g}_\sigma)$ is  viewed as a tangent vector to $ \mathcal R (M)$ at $\sigma$.

 In order to prove Theorem~\ref{convthmrepresentearth} we will use property (\ref{symplgradl}). More precisely, we will show
 that for any closed 1-form $d\xi$ with values in $ad(E_\sigma)$ 
\begin{eqnarray}\label{symplgradlnew}
(d\mbox{length}_{\lambda, \mu})_\sigma([d\xi])=dw(d\xi)=\Omega\left( [dw], [d\xi] \right).
\end{eqnarray}
Combined with (\ref{symplgradl}) and the fact that $\Omega$ is non-degenerate as a 2-form on $T_\sigma \mathcal R(M) \simeq H^1(M, ad(E_\sigma))$ it will imply the result. 
 We are going to prove this  for closed geodesics and then pass to the limit.
 
 The first step is to prove a formula for the first variation of the length function 
 \[
 \mbox{length}_\gamma: \mathcal R(M) \rightarrow \R
 \]
 for $\gamma$ a homotopy class of a simple closed curve in $M$. 
 
 \begin{lemma}\label{varleng} Let $\gamma$ be a homotopy class of a simple closed curve,  $\sigma: \pi_1(M) \rightarrow SO^+(2,1)$ a  representation defining a hyperbolic structure and  $\sigma(\gamma)=e^{lB}$ where $l$ is the length of the unique closed geodesic corresponding to $\gamma$ and $\sigma$. Let $[d\xi] \in T_\sigma \mathcal R(M)$ correspond to a  cocycle $\alpha$ via the De Rham Theorem~\ref{cohoisp}. Then
 \[
d(\mbox{length}_\gamma)_\sigma([d\xi])=1/2 \left(\alpha(\gamma), B \right)^\sharp.
 \]
 \end{lemma}
 \begin{proof} 
 By conjugating the representation if necessary we may assume, without loss of generality, that  $B \in \mathfrak s \mathfrak o(2,1)= \mathfrak s \mathfrak l(2,\R)$ is given by 
 $B=
\begin{pmatrix}
    & 1    & 0  \\                                
    & 0       & -1                                
\end{pmatrix}$.\
 Then 
\begin{equation}\label{lengthtra}
 Tr \sigma(\gamma) =2 \cosh l.
 \end{equation}

 Now let $\sigma_s=e^{s\alpha+O(s^2)}\sigma$ be a smooth family of representations  in $\mathcal R(M)$.  By differentiating (\ref{lengthtra}) with respect to $s$,
\begin{eqnarray*}
 2 l'(0)  \sinh l = Tr(\alpha(\gamma)e^{lB}).
 \end{eqnarray*}
Writing $\alpha(\gamma)=b B+[A,B]$ for a $A \in \mathfrak s \mathfrak o (2, 1) \simeq \mathfrak s \mathfrak l (2, \R)$, we verify
 \begin{eqnarray*}
 Tr(\alpha(\gamma)e^{lB})=  Tr(\alpha(\gamma)B) \sinh l.
 \end{eqnarray*} 
 By combining the above two equations
 \begin{equation*}
 l'(0)= 1/2 Tr(\alpha(\gamma)B)= 1/2 (\alpha(\gamma), B)^\sharp.
 \end{equation*} 
 \end{proof}
 \begin{remark}
 Our formulas may vary by a factor of 2 from other's.  This is due to the fact that our metric on $\HH$ is induced from the trace inner product and not 4 times the trace inner product which is the more standard convention.
 \end{remark}

  We next recall
   from \cite[Section 8.1]{daskal-uhlen2} the notion of parallel neighborhood around a closed geodesic $\gamma$.
   Parallel neighborhoods can be used instead of flow-boxes when the lamination consists of closed geodesics.  Given a closed geodesic $\gamma(t)=e^{tB}X_0$, we let $F: (0, 2\pi) \times (-\epsilon, \epsilon) \rightarrow \mathcal O=\mathcal O(\gamma) \subset M$
 so that the curves $s \mapsto  F (t,s)$ are geodesics parameterized by arc length, $F(t,0)=\gamma(t)$ and the  curves $t \mapsto F (t,s)$  are curves of fixed distance (parallel) to $\gamma$. 
Given  a smooth closed form $d\xi \in \Omega^1(M, Ad_\sigma(E))$,  we denote $\xi(t,s)=\xi(F(t,s))$ as a local function with values in the Lie algebra.

\begin{lemma}\label{formudvdphi}  Let $M$ be a closed surface with hyperbolic structure defined by a representation $\sigma: \pi_1(M) \rightarrow SO^+(2,1)$ and  $\gamma(t)=e^{tB}X_0$  a  closed geodesic in $M$. Let $dw=   \frac{d(\gamma(t))}{dt} \times \gamma(t)\delta(s)ds$ be a Lie algebra valued transverse measure on $\gamma$, where $\delta(s)$ is the delta-function supported on the geodesic.
 Let the cohomology class of $d\xi$ correspond to the cocycle $\alpha$ via the De Rham Theorem~\ref{cohoisp}. Then
\begin{eqnarray*}
dw(d\xi) =  (B, \alpha(\gamma))^\sharp.
\end{eqnarray*}
\end{lemma}
\begin{proof}
In the  parallel neighborhood $\mathcal O$,
\begin{eqnarray*}
dw(d\xi) &=&  \int_0^{2\pi} \left(B, \frac {d\xi(t,0)}{dt}\right)^\sharp dt \\
&=&  \int_0^{2\pi} \left(B, \xi(2\pi,0)-\xi(0,0)\right)^\sharp dt \\
&=&  \int_0^{2\pi} \left(B,  Ad(\sigma(\gamma))\xi(0,0)-\xi(0,0)+\alpha(\gamma)\right)^\sharp dt \\
&=& (B, \alpha(\gamma))^\sharp.
\end{eqnarray*}
The last equality follows by the $Ad$-invariance of the metric and the fact that, by definition of the representation $\sigma$, $\sigma(\gamma)=e^B$ commutes with $B$.  
\end{proof}

\begin{convthmrepresentearth}
Let  $(\gamma_k, \nu_k)$ be  a sequence of weighted closed geodesics converging in measure to $(\lambda, \mu)$. We have already shown that
\begin{eqnarray}\label{derilegt2}
d(\mbox{length}_{\gamma_k, \nu_k})_\sigma([d\xi]) =1/2dw_k ([d\xi])
\end{eqnarray}
and by \cite{ker1} and  \cite{ker2}
\begin{eqnarray}\label{dericonvrg}
d(\mbox{length}_{\gamma_k, \nu_k})_\sigma([d\xi])  \rightarrow (d\mbox{length}_{\lambda, \mu})_\sigma([d\xi]).
\end{eqnarray}
Moreover, by the convergence Theorem~\ref{convlemma}, $dw_k \rightharpoonup dw$.  Combined with  (\ref{derilegt2}) and (\ref{dericonvrg}) it implies (\ref{symplgradl}) for $(\lambda, \mu)$.
\end{convthmrepresentearth}


\begin{thebibliography}{ABC}

\bibitem[Ahl] {ahlfors} L. Ahlfors.  {\it Some remarks on Teichm\"uller's space of Riemann surfaces.} Ann. Math. 74 (1), 171-191 (1961).
\bibitem[Alb] {alberti}G. Alberti. {\it Rank one property for derivatives of functions with bounded variation}. Proc. Roy. Soc. Edinburgh Sect. A 123, no. 2, 239-274 (1993).
\bibitem[A-F-P]{ambrosio} L. Ambrosio, N. Fusco, D. Pallara. {\it Functions of bounded variation and free discontinuity problems.} Oxford Mathematical Monographs (2000).
\bibitem[Anz] {anzellotti} G. Anzellotti. {\it Pairings between measures and bounded functions and compensated compactness}. Annali di Matematica Pura ed Applicata 135.1, 293-318  (1983).
\bibitem[Ar1]{aronson1} G. Aronsson. {\it On certain singular solutions of the partial differential equation
$u_x^2u_{xx}+2u_xu_yu_{xy}+u_y^2u_{yy}=0$.}
Manuscripta math. 47, 133-151 (1984).
 \bibitem[Ar2]{aronson2} G. Aronsson. {\it Constructon of  singular solutions to the $p$-harmonic  equation and its limit equation for $p=\infty$}.
 Manuscripta math. 56, 135-158 (1986).
 \bibitem[Ar-C-J]{arcrju} G. Aronsson, M. Crandal and P. Juutinen. {\it A tour of the theory of absolutely minimizing functions.}
 Bulletin of the AMS., Vol. 41, No 4, 439-505 (2004).
  \bibitem[Ar-L]{aronsonlin} G. Aronsson and P. Lindqvist. {\it{On p-Harmonic Functions in the Plane and Their Stream Functions.}} Journal of Differential Equations, Vol. 74, Issue 1, 157-178 (1988).
\bibitem[Bac1]{aidan} A. Backus. {\it Minimal laminations and level sets of 1-harmonic functions.} Preprint.
  \bibitem[Bac2]{aidan2} A. Backus. {\it An infinity Laplacian for differential forms and calibrated laminations.} Preprint.
\bibitem[Bo1]{bonahon1} F. Bonahon. {\it Geodesic laminations on surfaces.} Contemporary Mathematics Volume 269 (2001).
\bibitem[Bo2]{bonahon3} F. Bonahon. {\it Shearing hyperbolic surfaces, bending pleated surfaces and Thurston's symplectic form.} Annales de la faculte des sciences de Toulouse, tome 5, no 2, 233-297 (1996).
\bibitem[C]{crandal} M. Crandall.
{\it A Visit with the $\infty$-Laplace Equation.} Calculus of Variations and Nonlinear Partial Differential Equations, 75-122 (2008).
 \bibitem[D-U1]{daskal-uhlen1} G. Daskalopoulos and K. Uhlenbeck.{ \it Transverse measures and  best Lipschitz and least gradient maps.} To appear J. Diff. Geom.
\bibitem[D-U2]{daskal-uhlen2} G. Daskalopoulos and K. Uhlenbeck.{ \it Analytic properties of minimal stretch maps and geodesic laminations.} Preprint.
\bibitem[D-D-G-S]{dacing} J. Danciger, T. Drumm, W. Goldman and I. Smilga. {\it Proper actions of discrete groups of affine transformations}. Dynamics, Geometry, Number Theory, University of Chicago Press, (2022).
\bibitem[Ek-Te]{temam}I. Ekeland and  R. Temam. {\it Convex Analysis and Variational Problems.} USA: Society for Industrial and Applied Mathematics (1999).
\bibitem[Go]{gold} W. Goldman. {\it The Symplectic Nature of Fundamental Groups of Surfaces.} Advances in Math.
Vol. 54, Issue 2,  200-225, (1984).
\bibitem[Gu-K]{kassel} F. Gueritaud and F. Kassel. {\it Maximally stretched laminations on geometrically finite hyperbolic manifolds.} Geom. Topol. Volume 21, Number 2, 693-840 (2017).
\bibitem [Ker1]{ker1} S. Kerckhoff. {\it The Nielsen realization problem.}  Annals of Mathematics 117, 235-265 (1983).
\bibitem [Ker2]{ker2} S. Kerckhoff. {\it Earthquakes are analytic. } Commentarii Mathematici Helvetici 60, 17-30, (1985).
\bibitem[Kir]{kirszbraun} M. Kirszbraun.  {\it \"Uber die zusammenziehende und Lipschitzsche Transformationen}. Fundamenta Mathematicae. 22: 77-108 (1934).
\bibitem[K-S]{korevaar-schoen1}  N. Korevaar and R. Schoen.  {\it Sobolev spaces and harmonic maps into metric space targets. }  Comm. Anal. Geom. 1  561-659, (1993).
\bibitem [Mc]{shane} E. McShane. {\it Extension of range of functions.} Bull. Amer. Math. Soc. 40,  837-842 (1934).
\bibitem[P-W]{wolf} H. Pan and  M. Wolf. {\it Ray structures on Teichm\"uller Space.} Preprint, arXiv: 2206.01371 (2022).
\bibitem[Pa-Th]{papa} A. Papadopoulos and G. Theret. {\it On Teichm\"uller's metric and Thurston's asymmetric metric on Teichm\"uller space.}  Handbook of Teichm\"uller Theory, Volume 1, 11, European Math. Soc. Publishing House (2007).
\bibitem[Ru-S]{sullivan} D. Ruelle and D. Sullivan. {\it Currents, flows and diffeomorphisms}. Topology 14, 319-327 (1975).
\bibitem[Si]{simon} L. Simon. {\it Introduction to Geometric Measure Theory.} Tsinghua Lectures (2014).
\bibitem[Thu1]{thurston} W. Thurston.
{\it Minimal stretch maps between hyperbolic surfaces.} Preprint arXiv:math/9801039.
\bibitem[Thu2]{thurston2} W. Thurston. {\it The Geometry and Topology of Three-Manifolds.} MSRI publications (2002).
\bibitem[Thu3]{thurston3} W. Thurston. {\it Earthquakes in 2-dimensional hyperbolic geometry.} In Fundamentals of hyperbolic geometry: selected expositions, volume 328, London Math. Soc. Lecture Note Ser., Cambridge Univ. Press, Cambridge, 2006.
\bibitem[Whi]{whitney} H. Whitney. {\it{Analytic extensions of differentiable functions defined in closed sets}}. Trans. Amer. Math. Soc. 36  no. 1, 63-89, (1934).
\end{thebibliography}
\end{document}